\documentclass[a4paper,12pt]{article}

\usepackage{amsmath,amssymb,amsfonts,graphicx}

\usepackage{epsfig}
\usepackage{color}
\usepackage{comment}
\usepackage{ulem}
\usepackage{mathdots}
\usepackage{graphics}
\usepackage{epstopdf}
\usepackage{mathdots}
\usepackage{epsfig}
\usepackage{verbatim}
\usepackage{color}
\usepackage{array}
\usepackage{url}

\def\IR{{\mathbb R}}
\def\IC{{\mathbb C}}

\def\IL{{\mathbb L}}

\newcommand{\sIL}{{{{\mathbb L}_s}}}

\def\sq{ { {\color{blue} \hfill\rule{1.5mm}{1.5mm} } } }

\newtheorem{lemma}{Lemma}[section]
\newtheorem{definition}{Definition}[section]
\newtheorem{proposition}{Proposition}[section]

\newtheorem{remark}{Remark}[section]

\newtheorem{example}{Example}[section]

\newcommand{\bA}{{\bf A}}
\newcommand{\bB}{{\bf B}}
\newcommand{\bC}{{\bf C}}
\newcommand{\bD}{{\bf D}}
\newcommand{\bE}{{\bf E}}

\newcommand{\bS}{{\bf S}}
\newcommand{\bY}{{\bf Y}}
\newcommand{\bJ}{{\bf J}}
\newcommand{\bK}{{\bf K}}

\newcommand{\bL}{{\bf L}}
\newcommand{\bM}{{\bf M}}
\newcommand{\bN}{{\bf N}}
\newcommand{\bI}{{\bf I}}
\newcommand{\bH}{{\bf H}}
\newcommand{\bP}{{\bf P}}
\newcommand{\bO}{{\bf O}}
\newcommand{\bQ}{{\bf Q}}
\newcommand{\bW}{{\bf W}}
\newcommand{\bR}{{\bf R}}
\newcommand{\bT}{{\bf T}}
\newcommand{\bZ}{{\bf Z}}
\newcommand{\bX}{{\bf X}}
\newcommand{\bx}{{\bf x}}
\newcommand{\by}{{\bf y}}
\newcommand{\bu}{{\bf u}}
\newcommand{\bV}{{\bf V}}
\newcommand{\bU}{{\bf U}}

\newcommand{\bfe}{{\bf e}}

\newcommand{\bn}{{\bf n}}

\newcommand{\bg}{{\bf g}}
\newcommand{\boo}{{\bf o}}

\newcommand{\br}{{\bf r}}

\newcommand{\bh}{{\bf h}}

\newcommand{\bv}{{\bf v}}
\newcommand{\bw}{{\bf w}}
\newcommand{\bz}{{\bf z}}

\newcommand{\bee}{{\bf e}}
\newcommand{\bff}{{\bf f}}
\newcommand{\bfz}{{\mathbf 0}}
\newcommand{\cC}{ {\cal C} }

\newcommand{\cO}{ {\cal O} }

\newcommand{\Si}{ \boldsymbol{\Sigma} }

\newcommand{\bPhi}{ \boldsymbol{\Phi} }

\newcommand{\bell}{\boldsymbol{\ell}}
\newcommand{\blambda}{\boldsymbol{\lambda}}

\newcommand{\bmu}{\boldsymbol{\mu}}

\newcommand{\bLambda}{\boldsymbol{\Lambda}}

\newcommand{\bXi}{\boldsymbol{\Xi}}

\newcommand{\cR}{ {\cal R} }

\newcommand{\cL}{ {\cal L} }

\title{
   Data-driven model order reduction of\\[.5mm] linear switched systems
}
\author{\\
I.V. Gosea \footnotemark[2]~, M. Petreczky\footnotemark[3]~,
   A.C. Antoulas\footnotemark[4]  \ \footnotemark[2] .
}

\date{}

\begin{document}

\maketitle

\renewcommand{\thefootnote}{\fnsymbol{footnote}}
\footnotetext[2]{
Data-Driven System Reduction and Identification Group, Max Planck Institute for Dynamics of Complex Technical Systems, Sandtorstrasse 1, 39106, Magdeburg, Germany
({\tt gosea@mpi-magdeburg.mpg.de})
}
\footnotetext[3]{Centre de Recherche en Informatique, Signal et Automatique de Lille (CRIStAL),  UMR CNRS 9189, CNRS, Ecole Centrale de Lille, France
({\tt mihaly.petreczky@ec-lille.fr})
}
\footnotetext[4]{
   Department of Electrical and Computer Engineering, 
   Rice University, 6100 Main St, MS-366, Houston, TX 77005, USA 
   ({\tt aca@rice.edu})
}

\renewcommand{\thefootnote}{\arabic{footnote}}

\begin{abstract} 
The Loewner framework for model reduction is extended to the class of linear switched systems. One advantage of this framework is that it introduces a trade-off between accuracy and complexity. Moreover, through this procedure, one can derive state-space models directly from data which is related to the input-output behavior of the original system. Hence, another advantage of the framework is that it does not require the initial system matrices. More exactly, the data used in this framework consists in frequency domain samples of input-output mappings of the original system. The definition of generalized transfer functions for linear switched systems resembles the one for bilinear systems. A key role is played by the coupling matrices, which ensure the transition from one active mode to another. 
\end{abstract}

\section{Introduction}
  
Model order reduction (MOR) seeks to transform large, complicated models of time dependent processes into smaller, simpler models that are nonetheless capable of representing accurately the behavior of the original  process under a variety of operating conditions. The goal is an efficient, methodical strategy that yields a dynamical system evolving in a substantially lower dimension space (hence requiring far less computational resources for realization), yet retaining  response characteristics close to the original system. Such reduced order models could be used as efficient surrogates for the original model, replacing it as a component in larger simulations.

Hybrid systems are a class of nonlinear systems which result from the interaction of continuous time dynamical sub-systems with discrete events. More precisely, a hybrid system is a collection of continuous time dynamical systems. The internal variable  of each dynamical system is governed by a set of differential equations. Each of the separate continuous time systems are labeled as a discrete mode. The transitions between the discrete states may result in a jump in the continuous internal variable. Linear switched systems (in short LSS) constitute a subclass of hybrid systems; the main property is that these systems switch among a finite number of linear subsystems. Also, the discrete events interacting with the sub-systems are governed by a piecewise continuous function called the switching signal.

Hybrid and switched systems are powerful models for distributed
embedded systems design where discrete controls are routinely
applied to continuous processes. However, the complexity
of verifying and assessing general properties of these systems
is very high so that the use of these models is limited in
applications where the size of the state space is large. To
cope with complexity, abstraction and reduction are useful
techniques. In this paper we analyze only the reduction part.

 In the past years, hybrid and switched systems have received increasing attention in the scientific community.
 For a detailed characterization of this relatively new class of dynamical systems, we refer the readers to the books \cite{book1LSS}, \cite{book2LSS}, \cite{book3LSS} and \cite{book4LSS}. Such systems are used in modeling, analysis and design of supervisory control systems, mechanical systems with impact, circuits with relays or ideal diodes.

The study of the properties of hybrid systems in general and switched systems in particular is still the subject of intense research, including the problems of stability (see \cite{dri02} and \cite{book2LSS}), realization including observability/controllability (see \cite{petre11} and \cite{ptt13}), analysis of switched DAE's (see \cite{mehrmann09} and \cite{tr12}) and numerical solutions (see \cite{mehrmann08}).
 Recently, considerable research has been dedicated to the problem of MOR for linear switched systems. The most prolific method that has been applied is balanced truncation (or some sort of gramian based derivation of it). Techniques that are based on balancing have been considered in the following: \cite{glw06}, \cite{ch09},\cite{bgmb10}, \cite{sw11}, \cite{mtc12}, \cite{pwl13} and \cite{pp16}. Also, another class of methods involve matching of generalized Markov parameters (known also as time domain Krylov methods) such as the ones in \cite{bastug16} and \cite{thesisB}; $\mathcal{H}_{\infty}$ type of reduction methods were developed in \cite{zsbw08}, \cite{bmb11} and \cite{zcw14}. Finally, we mention some publications that are focused on the reduction of discrete LSS, such as \cite{zbs09} and \cite{bmb12}.
 
  A linear switched system involves switching
between a number of linear systems (the modes of the
LSS). Hence, to apply balanced truncation
techniques to a switched linear system, one may seek for a
basis of the state space such that the corresponding
modes are all in balanced form.
It may happen that some state components of the LSS are difficult to reach and observe in some of the modes yet easy to reach and observe in others. In that case, deciding how to truncate the state variables and obtain a reduced
order model is not trivial. A solution to this problem is proposed
in \cite{mtc12} where it turns out that the average gramian can be used to obtain a reduced order model. This method will be used as a comparison tool for our new MOR method.

In the sequel we exclusively consider \textit {interpolatory} MOR  methods and in particular the Loewner framework applied to 
 LSS.
  Roughly speaking, in the linear case, interpolatory methods seek reduced models whose transfer function matches that of the original system at selected frequencies. For the nonlinear case, these methods require appropriate definitions of transfer functions.


In this paper, we focus on generalizing the Loewner Framework for  reducing linear switched systems. The presentation is tailored to emphasize the main procedure for a simplified case of LSS (i.e only two modes and LTI's in SISO format activating in both modes). The paper is organized as follows. In the next section, we review the formal definition of continuous-time linear switched systems. Furthermore, we introduce the generalized transfer functions for LSS as input-output mappings in frequency domain.
 Section 3 includes a brief introduction of the Loewner framework for linear systems. In
Section 4, we introduce the Loewner framework for LSS with two modes. In section 5, we generalize most of the results in the previous section for the case of LSS with $D \geqslant 2$ modes. Finally, in Section 6, we discuss the applicability  of the new introduced method for reducing LSS. In this sense, by means of three numerical examples (one of which large scale), we compare the results obtained by applying the Loewner method against the method in \cite{mtc12}. In Section 7, we present a summary of the findings and the conclusions.

\section{Linear switched systems} \label{sec:switched}

\begin{definition}
 A continuous time linear switched system (LSS) is a control system of the form:
\begin{equation}\label{LSS_def}
\Si: \begin{cases}
 \bE_{\sigma(t)} \dot{\bx}(t) = \bA_{\sigma(t)} \bx(t) + \bB_{\sigma(t)}u(t), \ \ \bx(t) = \bx_0, \\
\by(t) = \bC_{\sigma(t)} \bx(t),
\end{cases}
\end{equation}
where $Q = \{1,2,\ldots,D\}, D > 1$, is a set of discrete modes, $\sigma(t)$ is  the switching signal, $u$ is  the input, $\bx$ is the state, and y  is the output.

The system matrices  $\bE_q, \bA_q \in \mathbb{R}^{n_q \times n_q}, \ \bB_q \in \mathbb{R}^{n_q \times m}, \ \bC_q \in \mathbb{R}^{p \times n_q}$, where $q \in Q$, correspond to the linear system active in mode $q \in Q$, and $\bx_0$ is the initial state. We consider the $\bE_q$ matrices to be invertible. Furthermore, the transition from one mode to another is made via the so called switching or coupling matrices $\bK_{q_1,q_2} \in \mathbb{R}^{n_{q_2} \times n_{q_1}}$ where  $q_1,q_2 \in Q$. 
\end{definition}

\begin{remark}
{\rm The case for which the coupling is made between identical modes is excluded, Hence, when $q_1 = q_2 = q$, consider that the coupling matrices are identity matrices, i.e. $\bK_{q,q} = \bI_{n_q}$.}
\end{remark}

The notation $\Si = (n_1,n_2,\ldots,n_D, \{(\bE_q,\bA_q, \bB_q, \bC_q)| q \in Q\},\{ \bK_{q_{i},q_{i+1}} | q_i, q_{i+1} \in Q \}, \bx_0)$ is used as a short-hand representation for LSS's described by the equations in (\ref{LSS_def}). The vector $\bn = \left( \begin{array}{cccc}
n_1 & n_2 & \cdots & n_D
\end{array} \right)$ is the dimension (order) of $\Si$. The linear system which is active in the $q^{th}$ mode of $\Si$ is denoted with $\Si_q$ and it is described by (where $1 \leqslant q \leqslant D$)
\begin{equation}\label{LSS_def2}
\Si_k: \begin{cases}
\bE_{q} \dot{\bx_q}(t) = \bA_{q}\bx_q(t) + \bB_{q}u(t), \ \ \bx(t_k) = \bx_k,  \\
y(t) = \bC_{q} \bx_q(t).
\end{cases}
\vspace{-2mm}
\end{equation}

The restriction of the switching signal $\sigma(t)$ to a finite interval of time $[0, T]$ can be interpreted as finite sequence of elements of $Q \times \mathbb{R}_{+}$ of the form:
$$
\nu(\sigma) = (q_1,t_1) (q_2,t_2) \ldots (q_k,t_k),
$$
where $q_1, \ldots , q_k \in Q$ and $0<t_1 < t_2 < \cdots < t_k  \in \mathbb{R}_{+}, \ t_1 + \cdots + t_k = T$, such that for all $ t \in [0, T]$ we have:
$$
\sigma(t) = \left\{ \begin{array}{l} q_1 \ \  \text{if} \ \ t \in [0,t_1], \\
q_2 \ \  \text{if} \ \  t \in (t_1,t_1+t_2], \\ 
\vdots \\
q_i \ \  \text{if} \ \ t \in (t_1+\ldots+t_{i-1},t_1+\ldots+t_{i-1}+t_i], \\
\vdots \\
q_k \ \  \text{if} \ \ t \in (t_1+\ldots+t_{k-1},t_1+\ldots+t_{k-1}+t_k]. \end{array} \right.
$$
In short, by denoting $T_i := t_1+\ldots+t_{i-1}+t_i, \  T_0 :=0, \  T_k := T$, write:
$$
\sigma(t) = \begin{cases}  q_1 \ \  \text{if} \ \ t \in [0,T_1], \\
 q_i \ \  \text{if} \ \ t \in (T_{i-1},T_i], \ i > 2. \end{cases}
$$

Denote by $PC(\mathbb{R}_{+},\mathbb{R}^n)$, $P_c(\mathbb{R}_{+},\mathbb{R}^n)$, the set of all piecewise-continuous, and piecewise-constant functions, respectively.
\begin{definition}
A tuple $(\bx,\bu,\sigma,\by)$, where $\bx:\mathbb{R}_{+} \rightarrow \bigcup_{i=1}^{D} \mathbb{R}^{n_i}$, $\bu \in PC(\mathcal{R}_{+},\mathbb{R}^{m}), \sigma \in P_c(\mathbb{R}_{+},Q), \by \in PC(\mathbb{R}_{+},\mathbb{R}^{p})$ is called a solution, if the following conditions simultaneously hold:

\begin{enumerate}
\item The restriction of $\bx(t)$  to $(T_{i-1},T_{i}]$ is differentiable, and satisfies $\bE_{q_i} \dot \bx (t) =\bA_{q_i}\bx(t)+\bB \bu(t)$.
\item  Furthermore, when switching from mode $q_{i}$ to mode $q_{i+1}$ at time $T_i$, the following holds
\vspace{-2mm}
\begin{equation*}
\bE_{q_{i+1}} \displaystyle \lim_{t \searrow T_i} \bx_{q_{i+1}}(t) = \bK_{q_i,q_{i+1}}  \bx_{q_i}(T_i).
\vspace{-2mm}
\end{equation*}  
\item  Moreover, for all $t \in \mathbb{R}$, $\by(t)=\bC_{\sigma(t)} \bx(t)$ holds.
\end{enumerate}
\end{definition}

\vspace{-2mm}

The switching matrices $\bK_{q_{i},q_{i+1}}$ allow having different dimensions for the subsystems active in different modes. For instance, the pencil $(\bA_{q_i},\bE_{q_i}) \in \mathbb{R}^{n_{q_i} \times n_{q_i}}$, while the pencil $(\bA_{q_{i+1}},\bE_{q_{i+1}}) \in \mathbb{R}^{n_{q_{i+1}} \times n_{q_{i+1}}}$ where the values $n_{q_i}$ and $n_{q_{i+1}}$ need not be the same. 
If the $\bK_{q_{i},q_{i+1}}$ matrices are not explicitly given, it is considered that they are identity matrices.

The input-output behavior of an LSS system can be formalized in time domain as a map $\bff(\bu,\sigma)(t)$. This particular map can be written in \textit{generalized kernel representation} (as suggested in \cite{pvs10}) using the unique family of analytic functions: $\bg_{q_1,\ldots,q_k}: \mathbb{R}_{+}^k \rightarrow \mathbb{R}^p$ and $\bh_{q_1,\ldots,q_k}: \mathbb{R}_{+}^k \rightarrow \mathbb{R}^{p \times m}$ with $q_1,\ldots,q_k \in Q, \ k \geqslant 1$ such that for all pairs $(\bu,\sigma)$ and for $T = t_1+t_2+\cdots+t_k$ we can write:
$$
\bff(u,\sigma)(t) = \bg_{q_1,q_2,\ldots,q_k}(t_1,t_2,...,t_k)+ \displaystyle \sum_{i=1}^{k} \int_{0}^{t_i} \bh_{q_i,q_{i+1},\ldots,q_k}(t_i-\tau,t_{i+1},\ldots,t_k) u(\tau+T_{i-1}) d\tau,
$$
where the functions $\bg, \bh$ are defined for $k \geqslant 1$, as follows,
\begin{equation}\label{init_state}
\bg_{q_1,q_2,\ldots,q_k}(t_1,t_2,\ldots,t_k) = \bC_{q_k}e^{\tilde{\bA}_{q_k}t_k} \tilde{\bK}_{q_{k-1},q_{k}} e^{\tilde{\bA}_{q_{k-1}}t_{k-1}} \tilde{\bK}_{q_{k-2},q_{k-1}} \cdots \tilde{\bK}_{q_{1},q_{2}} e^{\tilde{\bA}_{q_1}t_1} \bx_0,
\end{equation}
\begin{equation}\label{kern}
\bh_{q_1,q_2,\ldots,q_k}(t_1,t_2,\ldots,t_k) = \bC_{q_k}e^{\tilde{\bA}_{q_k}t_k} \tilde{\bK}_{q_{k-1},q_{k}} e^{\tilde{\bA}_{q_{k-1}}t_{k-1}} \tilde{\bK}_{q_{k-2},q_{k-1}} \cdots \tilde{\bK}_{q_{1},q_{2}} e^{\tilde{\bA}_{q_1}t_1}\tilde{\bB}_1.
\end{equation}
Note that, for the functions defined in (\ref{init_state}) and (\ref{kern}) we consider the $\bE_{q_i}$ matrices to be incorporated into the $\bA_{q_i}$ and $\bB_{q_i}$ matrices (i.e. $\tilde{\bA}_{q_i} = \bE_{q_i}^{-1} \bA_{q_i}, \ \tilde{\bB}_{q_i} = \bE_{q_i}^{-1} \bB_{q_i}$). Moreover, the transformed coupling matrices are written accordingly $\tilde{\bK}_{q_i,q_{i+1}} = \bE_{q_{i+1}}^{-1} \bK_{q_i,q_{i+1}}$. \\

\vspace{-2mm}

In the rest of the paper, the LSS we treat {\it are assumed to have zero initial conditions}, i.e., $\bx_0 = \bfz$.  Hence, only the $\bh$ functions in (\ref{kern}) are relevant for characterizing the input-output mapping $\bff$.

The behavior of the input-output mappings in frequency domain is in turn characterized by a series of multivariate rational functions obtained by taking the multivariable Laplace transform of the regular kernels in (\ref{kern}), as for


$$
\bH_{q_1}(s_1) = \bC_{q_1} \bPhi_{q_1}(s_1) \bB_{q_1},
$$
$$
\bH_{q_1,q_2}(s_1,s_2) = \bC_{q_1} \bPhi_{q_1}(s_1) \bK_{q_2,q_1}\bPhi_{q_2}(s_2)\bB_{q_2},
$$
$$
\bH_{q_1,q_2,q_3}(s_1,s_2,s_3) = \bC_{q_1} \bPhi_{q_1}(s_1) \bK_{q_2,q_1}\bPhi_{q_2}(s_2) \bK_{q_3,q_2} \bPhi_{q_3}(s_3) \bB_{q_3}, \ \  \cdots
$$
\normalsize
For $k \geqslant 1$, write the {\it level k generalized transfer function} associated to the switching sequence $(q_1,q_2,\ldots,q_k)$, and evaluated at the points $(s_1,s_2,\ldots s_k)$ as,
\begin{equation}\label{gen_trf_lss}
\bH_{q_1,q_2,...,q_k}(s_1,s_2,...,s_k) = \bC_{q_1} \bPhi_{q_1}(s_1) \bK_{q_2,q_1} \bPhi_{q_2}(s_2) \cdots \bK_{q_{k},q_{k-1}}\bPhi_{q_k}(s_k) \bB_{q_k},
\end{equation}
where $\bPhi_q(s) = (s \bE_q - \bA_q)^{-1}$, $q_j \in \{1,2,...,D\}, \ 1 \leqslant j \leqslant k$ and $k \geqslant 3$. These functions are the {\it generalized transfer functions} of the linear switched system $\Si$. Their definition is similar to the ones corresponding to bilinear systems (see \cite{agi16}).
  
By using their samples, we are able to directly come up with (reduced) switched models that interpolate the original model - generalization of the Loewner framework to LSS. 

We construct LSS reduced models by means of matching samples of input-output mappings corresponding to the original LSS system and evaluated at finite sampling points (as opposed to other approaches - see \cite{thesisB} and \cite{bastug16}, where the behavior at infinity is studied instead, i.e. by matching Markov parameters).

 For the explicit derivation of these types of transfer functions (which is based on the so-called Volterra series representation) we refer the readers to \cite{Ru82}.\\

\section{Interpolatory MOR methods and the Loewner framework} 

 
Consider a full-order linear system defined by  $\bE \in \mathbb{R}^{n \times n}$, $\bA \in \mathbb{R}^{n \times n}$, $\bB \in \mathbb{R}^{n \times m}$, $\bC \in \mathbb{R}^{p \times n}$, and its transfer function $\bH(s)$ $=$ $\bC(s\bE-\bA)^{-1}\bB$. Given {\it  left interpolation points}:  $\{\mu_j\}_{j=1}^q$ $\subset$ $\IC$, with {\it  left tangential directions}:  $\{{\bell}_j\}_{j=1}^q$ $\subset$ $\IC^p$, and {\it right interpolation points}: $\{\lambda_i\}_{i=1}^k$ $\subset$ $\IC$, with {\it right tangential directions}: $\{{\br}_i\}_{i=1}^k$ $\subset$ $\IC^m$, find a reduced-order system
$\hat\bE, \hat\bA, \hat\bB, \hat\bC$, 
such that the resulting transfer function, $\hat\bH (s)$ is a {\it tangential interpolant} to $\bH(s)$:
\vspace{-1mm}
\begin{equation} \label{prob1}
\begin{array}{ccc}
{\bell}_j^T \hat\bH(\mu_j)  =  {\bell}_j^T \bH(\mu_j),~ j=1,\ldots,q,
 \ \ \mbox{and} \ \ \hat\bH(\lambda_i) {\br}_i =  \bH(\lambda_i) {\br}_i,~ i=1,\ldots,k
\end{array}
\vspace{-2mm}
\end{equation}
Interpolation points and tangent directions are selected to realize appropriate MOR goals. If instead of state space data, we are given {\it input/output data}, the resulting problem is hence modified. Given a set of input-output response measurements specified by {\it  left driving frequencies}: $\{\mu_j\}_{j=1}^q$ $\subset$ $\IC$, using {\it left input directions}: $\{{\bell}_j\}_{j=1}^q$ $\subset$ $\IC^p$, producing {\it left responses}:  $\{{\bv}_j\}_{j=1}^q$ $\subset$ $\IC^m$, and {\it right driving frequencies}: $\{\lambda_i\}_{i=1}^k$ $\subset$ $\IC$, using {\it right input directions}:  $\{{\br}_i\}_{i=1}^k$ $\subset$ $\IC^m$, producing {\it right responses}:  $\{{\bw}_i\}_{i=1}^k$ $\subset$ $\IC^p$, find (low order) system matrices $\hat\bE$, $\hat\bA$, $\hat\bB$, $\hat\bC$, such that the resulting transfer function, $\hat\bH (s)$, is an ({\it approximate}) {\it tangential interpolant} to the data:
\vspace{-1mm}
\begin{equation} \label{prob2} 
{\bell}_j^T \hat\bH(\mu_j) ~ = 
~{\bv}_j^T, \ \ \ j=1,\ldots,q, ~~~\mbox{and}~~~ 
\hat\bH(\lambda_i) {\br}_i~ = 
~{\bw}_i, \ \ \ i=1,\ldots,k.
\vspace{-1mm}
\end{equation}



\subsection{Overview of the Loewner framework for linear systems}

The approach we discuss in this section is {\it data driven}. After collecting input/output (e.g. frequency response) measurements for some appropriate range of frequencies, we construct models which fit (or approximately fit) the data and have reduced dimension.
{\it The key is that,
larger amounts of data than necessary are collected
and the essential underlying system structure is extracted appropriately.}
Thus an {\it advantage} of this approach is that it can provide the user with a trade-off between accuracy of fit and complexity of the model.

 The {\it Loewner framework} was developed in a series of papers;
for details we refer the reader to \cite{book}, as well as \cite{ajm07, left_aca_11, tcad09, ail12, IonAnt13, IonAnt14}. For a recent overview see \cite{birkjour}.

\subsection{The Loewner pencil}

We will formulate the results for the more general {\it tangential interpolation} problem. 
We are given the {\it right data}: ~$(\lambda_i;\br_i,\bw_i)$, $i=1,\cdots,k$, ~ and the
{\it left data}: ~$(\mu_j;\bell_j^T,\bv_j^T)$, $j=1,\cdots,q$;~ it is assumed 
for simplicity that all points are distinct. The dimensions are as in (\ref{prob1}), (\ref{prob2}).
The data can be organized as follows: the {\it right data}:
$$
\bLambda=\mbox{diag}\left[
\lambda_1, \ldots,\lambda_k \right]\in\IC^{k\times k},~
\bR=[\br_1,\ldots,\br_k] \in\IC^{m\times k},~
\bW=[\bw_1,\ldots,\bw_k]\in\IC^{p\times k},
$$
and the {\it left data}:
$$
\bM=\mbox{diag}\left[ \mu_1,\ldots, \mu_q \right] \in \IC^{q\times q} ,~
\bL^T=\left[ \bell_1,\ldots,\bell_q\right] \in \IC^{q\times p} ,~
\bV^T=\left[ \bv_1,\ldots,\bv_q\right]\in\IC^{q\times m} .
$$
Then, the associated {\it Loewner pencil}, consists of the {\it Loewner and shifted Loewner matrices}.
The {\bf Loewner matrix ~$\IL \in\IC^{q\times k}$},~ is defined as:
$$
\IL=\left[\begin{array}{ccc}
\frac{\bv_1^T\br_1-\bell_1^T\bw_1}{\mu_1-\lambda_1} & \cdots &
\frac{\bv_1^T\br_k-\bell_1^T\bw_k}{\mu_1-\lambda_k} \\
\vdots & \ddots & \vdots \\
\frac{\bv_q^T\br_1-\bell_q^T\bw_1}{\mu_q-\lambda_1} & \cdots &
\frac{\bv_q^T\br_k-\bell_q^T\bw_k}{\mu_q-\lambda_k} \\
\end{array}\right]
$$
$\IL$ satisfies the Sylvester equation ~$
\IL\bLambda-\bM \IL= \bV\bR-\bL\bW$.~ 
Suppose that the underlying transfer function is ~$\bH(s)=\bC(s\bE-\bA)^{-1}\bB$,~ and define
the {\it generalized observability/controllability} matrices:
\small
\begin{equation} \label{eq:OR_linear}
\cO=
\left[\!\begin{array}{c}
\bC(\mu_1\bE-\bA)^{-1}\\
\vdots\\
\bC(\mu_q\bE-\bA)^{-1}
\end{array}
\!\right],~~~\cR=
\left[\!
\begin{array}{ccc}
(\lambda_1\bE-\bA)^{-1}\bB & \cdots & (\lambda_k\bE-\bA)^{-1}\bB
\end{array}\!\right].
\end{equation}
\normalsize
It readily follows that the Loewner matrix can be factored as
$\IL=-\cO\bE\cR$. The {\bf shifted Loewner matrix ~$\sIL \in\IC^{q\times k}$},~ is defined as:
$$
\sIL=\left[\begin{array}{ccc}
\frac{{\mu_1}\bv_1^T\br_1-\bell_1^T\bw_1{\lambda_1}}{\mu_1-\lambda_1} & \cdots &
\frac{{\mu_1}\bv_1^T\br_k-\bell_1^T\bw_k{\lambda_k}}{\mu_1-\lambda_k} \\
\vdots & \ddots & \vdots \\
\frac{{\mu_q}\bv_q^T\br_1-\bell_q^T \bw_1{\lambda_1}}{\mu_q-\lambda_1} & \cdots &
\frac{{\mu_q}\bv_q^T\br_k-\bell_q^T\bw_k{\lambda_k}}{\mu_q-\lambda_k} \\
\end{array}\right] 
$$
$\sIL$ satisfies the Sylvester equation \,
$
\sIL\bLambda-\bM \sIL= \bM\bV\bR-\bL\bW \bLambda 
$,~ 
and can be factored in terms of the generalized controllability/observabilty matrices as $\sIL=-\cO\bA\cR$.
Finally notice that the following relations hold: $\bV=\bC\cR,~~ \bW=\cO\bB$.


\subsection{Construction of reduced order models}

We will distinguish two cases namely, the {\it right amount of data} and the more realistic {\it redundant amount of data} cases. The following lemma covers the first case.

\vspace{0mm}

\begin{lemma}{Assume that $k=q$, and let $\,(\sIL,\,\IL)$,\, be a regular pencil, such that none of the interpolation points $\lambda_i$, $\mu_j$ are its eigenvalues. Then $\bE=-\IL$, ~$\bA=-\sIL$,~ $\bB=\bV$, ~$\bC=\bW$,~ is a minimal realization of an interpolant of the data, i.e., the rational function ~$\bH(s)=\bW(\sIL-s\IL)^{-1}\bV$,~ interpolates the data (the conditions in (\ref{prob2}) are hence matched).}
\end{lemma}

\vspace{0mm}

If the pencil $(\sIL,\,\IL)$ is singular we are dealing with the case of {\it redundant data}. In this case if the following assumption is satisfied:
\vspace{-1mm}
\begin{equation} \label{condition}
\mbox{rank}\,(x\IL-\sIL) = \mbox{rank}\,\left[\begin{array}{l}{\IL}\\{\sIL}\end{array}\right]
=\mbox{rank}\,\left[{\IL}~~{\sIL}\right]=  r \leqslant k,
\vspace{-1mm}
\end{equation}
for all $x\in\{\lambda_i\}\cup\{\mu_j\}$,~ we consider the following SVD factorizations: 
\vspace{-1mm}
\begin{equation} \label{projection}
\left[\IL ~~~\sIL\right]=\bY_1 \bS_1\bX_1^T,~~
\left[\begin{array}{l}\IL\\\sIL\end{array}\right] = \bY_2 \bS_2 \bX_2^T ,
\vspace{-1mm}
\end{equation}
where ~$ \bY_1$, $\bX_2$ $\in$ $\IC^{k \times k}$. By selecting the first $r$ columns of the matrices $\bY_1$ and $\bX_2$, we come up with projection matrices $\bY,\bX \in \mathbb{C}^{k \times r}$. The following result is used in practical applications.\\
\vspace{-4mm}
\begin{lemma}
A realization $(\bE,\bA,\bB,\bC)$ of an approximate interpolant is given by the system matrices $\bE = -\bY^T \IL \bX$,~ $\bA = -\bY^T\sIL \bX$, ~$\bB = \bY^T\bV$,~ $\bC = \bW \bX$. Hence, the rational function ~$\bH(s)=\bW \bX(\bY^T \sIL \bX-s \bY^T \IL \bX)^{-1}\bY^T \bV$ approximately matches the data (the conditions in (\ref{prob2}) are approximately fulfilled, i.e. $\bH(\lambda_i) \br_i = \bw_i+\epsilon^r_i$ and $\bell_j^T \bH(\mu_j) = \bv_j^T+(\epsilon^\ell_j)^T$, where the residual errors are collected in the vectors $\epsilon^r_i$ and $\epsilon^\ell_j$).
\end{lemma}

\vspace{-1mm}

Thus, if we have more data than necessary, we can consider ~$(\sIL,\,\IL,\,\bV,\,\bW)$,~ as a {\it singular} model of the data. An appropriate projection then yields a reduced
system of order $k$ (see \cite{ajm07,drazin}).

A direct consequence is that the {\it Loewner framework offers a trade-off between accuracy and complexity of the reduced order system}, by means of the singular values of $\IL$.

\begin{remark}
{\rm For an error bound that links the quality of approximation to the singular values of the Loewner pencil (which is valid only at the interpolation points $\mu_j$ and $\lambda_i$), we refer the readers to \cite{birkjour}.}
\end{remark}

\section{The Loewner framework for LSS - the case D=2}

The characterization of linear switched systems by means of rational functions suggests that reduction of such 
systems can be performed by means of {\it interpolatory methods}. In the following we will show how to generalize the  Loewner framework to LSS by interpolating appropriately defined transfer functions on a chosen grid of frequencies (interpolation points). 

As for the linear case, the given set of sampling (interpolation) points is first partitioned into the two following categories: {\it  left interpolation points}:  $\{\mu_j\}_{j=1}^{\ell}$ $\subset$ $\IC$ and and {\it right interpolation points}: $\{\lambda_i\}_{i=1}^k$ $\subset$ $\IC$. In this paper we consider only the case of SISO linear switched systems- hence the  {\it  left} and {\it right tangential directions} can be considered to be scalar (i.e. taking the value 1). Since the transfer functions which are going to be matched are not single variable functions anymore (they depend on multiple variables as described in ), the structure of the interpolation points used in the new framework is going to change. Instead of having singleton values as in Section 3, we will use instead n-tuples that include multiple singleton values.

For simplicity of the exposition, we first consider the simplified case $D = 2$ (the LSS system switches between two modes only). This situation is encountered in most of the numerical examples in the literature we came across.
 Nevertheless, all the results presented in this section can be generalized for higher number of modes in a more or less straightforward way (as presented in Section 5).
Depending on the switching signal $\sigma(t)$, we either have,
$$
\Si_1: \ \begin{cases}
 \bE_{1} \dot{\bx}_1(t) = \bA_{1}\bx_1(t) + \bB_{1}u(t), \\
y(t) = \bC_{1} \bx_1(t) 
\end{cases} \ \text{or} \ \ \ \Si_2: \ \begin{cases}
\bE_{2} \dot{\bx}_2(t) = \bA_{2}x_2(t) + \bB_{2}u(t), \\
y(t) = \bC_{2} \bx_2(t) 
\end{cases}
$$
where $\text{dim}(\Si_1) = n_1$ (i.e. $\bx_1 \in \mathbb{R}^{n_1}$ and $\bE_1,\bA_1 \in \mathbb{R}^{n_1 \times n_1}, \bB_1,\bC_1^T \in \mathbb{R}^{n_1}$) and also $\text{dim}(\Si_2) = n_2$ (i.e. $\bx_2 \in \mathbb{R}^{n_2}$ and $\bE_2,\bA_2 \in \mathbb{R}^{n_2 \times n_2}, \bB_2,\bC_2^T \in \mathbb{R}^{n_2}$). Notice that we allow both the two subsystems to be written in descriptor format (having possibly singular E matrix).

Denote, for simplicity, with $\bK_1$ the coupling matrix when switching from mode 1 to mode 2 (instead of $\bK_{1,2}$) and, with $\bK_2$, the coupling matrix when switching from mode 2 to mode 1 (instead of $\bK_{2,1}$) with $\bK_1 \in \mathbb{R}^{n_2 \times n_1}$ and $\bK_2 \in \mathbb{R}^{n_1 \times n_2}$. 


The generalized transfer functions are defined as (where $\bPhi_q(s) = (s \bE_q-\bA_q)^{-1}, \ q \in \{1,2\}$),

\vspace{5mm}

\noindent
\underline{\textbf{Level 1}}\\

$ \bH_1(s_1) = \bC_1 \bPhi_1(s_1)\bB_1$ \ \ \ \ \ \ \ \ $ \bH_2(s_2) = \bC_2 \bPhi_2(s_2)\bB_2$\\

\noindent
\underline{\textbf{Level 2}}\\

$ \bH_{1,2}(s_1,s_2) = \bC_1 \bPhi_1(s_1) \bK_2 \bPhi_2(s_2)  \bB_2$ \ \ \ \ \ \ \ \ $ \bH_{2,1}(s_2,s_1) = \bC_2 \bPhi_2(s_2) \bK_1 \bPhi_1(s_1)  \bB_1$\\

\noindent
\underline{\textbf{Level 3}}

$$
 \begin{cases}
 \bH_{1,2,1}(s_1,s_2,s_3) = \bC_1 \bPhi_1(s_1) \bK_2 \bPhi_2(s_2) \bK_1 \bPhi_1(s_3)  \bB_1 \\
 \bH_{2,1,2}(s_1,s_2,s_3) = \bC_2 \bPhi_2(s_1) \bK_1 \bPhi_1(s_2) \bK_2 \bPhi_2(s_3)  \bB_2, \ \ \cdots
\end{cases} 
$$

%
%

\begin{definition}

Consider the two LSS, $\hat{\Si} = (n_1,n_2, \{(\hat{\bE}_i,\hat{\bA}_i, \hat{\bB}_i, \hat{\bC}_i)\}_{i=1}^2,\{ \hat{\bK}_{i,j}  \}_{i,j =1}^2, \bfz)$ and $\bar{\Si} = (n_1,n_2, \{(\bar{\bE}_i,\bar{\bA}_i, \bar{\bB}_i, \bar{\bC}_i)\}_{i=1}^2,\{ \bK_{i,j}  \}_{i,j =1}^2, \bfz)$. These systems are said to be equivalent if there exist non-singular matrices $\bZ_j^L$ and $\bZ_j^R$ so that
$$
\bar{\bE}_j = \bZ_j^L \hat{\bE}_j \bZ_j^R, \ \  \bar{\bA}_j = \bZ_j^L \hat{\bA}_j \bZ_j^R, \ \ \bar{\bB}_j = \bZ_j^L \hat{\bB}_j , \ \ \bar{\bC}_j = \hat{\bC}_j \bZ_j^R, \ \ j \in \{1,2\}
$$ 
and also $\bar{\bK}_1 = \bZ_2^L \hat{\bK}_1 \bZ_1^R, \ \bar{\bK}_2 = \bZ_1^L \hat{\bK}_2 \bZ_2^R$. In this configuration, one can easily show that the transfer functions defined above are the same for each LSS and for all sampling points $s_k$.
\end{definition}

\subsection{The generalized controllability and observability matrices}

Consider a LSS system $\Si$ as described in (\ref{LSS_def}) with $\text{dim}(\Si_k) = n_k$ for $k = 1,2$ and let $\bK_1 \in \mathbb{R}^{n_2 \times n_1}$ and $\bK_2 \in \mathbb{R}^{n_1 \times n_2}$ be the coupling matrices. Before stating the general definitions, we first clarify how the newly introduced matrices are constructed through a simple self-explanatory example.
\begin{example}
Consider 10 left interpolation points $
\{\mu_1^{(1)},\ldots,\mu_4^{(1)},\mu_1^{(2)},\ldots,\mu_6^{(2)}\}
$ which are written in nested multi-tuple format (corresponding to each mode of the LSS):
\begin{equation*}
\underline{\textbf{\text{Mode} 1}} : \ \ 
\begin{array}{ll}
 \!\bmu_1^{(1)}\!=\!\left\{\!\! \begin{array}{l}
\big{(} \mu^{(1)}_{1} \big{)},\,\\[2mm]
\big{(} \mu^{(1)}_2,\,\mu^{(1)}_{3} \big{)}
\end{array} \right.
&
 \!\bmu_1^{(2)}\!=\!\left\{\!\! \begin{array}{l}
\big{(} \mu^{(2)}_{1} \big{)},\,\\[2mm]
\big{(} \mu^{(2)}_2,\,\mu^{(2)}_{3} \big{)},\,\\[2mm]
\big{(} \mu^{(2)}_1,\,\mu^{(2)}_{4},\mu^{(2)}_{5} \big{)}
\end{array} \right.
\end{array}\!\!\!\!\! 
\end{equation*}
\begin{equation*}
\underline{\textbf{\text{Mode} 2}} : \ \ \begin{array}{ll}
 \!\bmu_2^{(1)}\!=\!\left\{\!\! \begin{array}{l}
\big{(} \mu^{(1)}_{2} \big{)},\,\\[2mm]
\big{(} \mu^{(1)}_1,\,\mu^{(1)}_{4} \big{)},\
\end{array} \right.
&
 \!\bmu_2^{(2)}\!=\!\left\{\!\! \begin{array}{l}
\big{(} \mu^{(2)}_{2} \big{)},\,\\[2mm]
\big{(} \mu^{(2)}_1,\,\mu^{(2)}_{4} \big{)},\,\\[2mm]
\big{(} \mu^{(2)}_2,\,\mu^{(2)}_{3},\mu^{(2)}_{6} \big{)}\,
\end{array} \right.
\end{array}\!\!\!\!\!
\end{equation*}
We explicitly write the generalized observability matrices $\mathcal{O}_1$ and $\mathcal{O}_2$ as follows:
$$
\cO_1= \left[ \begin{array}{c}  \bC_1\,\bPhi_1({\mu^{(1)}_{1}}) \\ \bC_2\,\bPhi_2({\mu^{(1)}_{2}})\,\bK_1\,\bPhi_1({\mu^{(1)}_{3}}) \\ \bC_1\,\bPhi_1({\mu^{(2)}_{1}}) \\ \bC_2\,\bPhi_2({\mu^{(2)}_{2}})\,\bK_1\,\bPhi_1({\mu^{(2)}_{3}}) \\ \bC_1\,\bPhi_1({\mu^{(2)}_{1}})\,\bK_2\,\bPhi_2({\mu^{(2)}_{4}}) \,\bK_1\,\bPhi_1({\mu^{(2)}_{5}})
 \end{array} \right], \ \ \ \cO_2= \left[ \begin{array}{c}  \bC_2\,\bPhi_2({\mu^{(1)}_{2}}) \\ \bC_1\,\bPhi_1({\mu^{(1)}_{1}})\,\bK_2\,\bPhi_2({\mu^{(1)}_{4}}) \\ \bC_2\,\bPhi_2({\mu^{(2)}_{2}}) \\ \bC_1\,\bPhi_1({\mu^{(2)}_{1}})\,\bK_2\,\bPhi_2({\mu^{(2)}_{4}}) \\ \bC_2\,\bPhi_2({\mu^{(2)}_{2}})\,\bK_1\,\bPhi_1({\mu^{(2)}_{3}}) \,\bK_2\,\bPhi_2({\mu^{(2)}_{6}})
 \end{array} \right]
$$
\end{example}

\begin{definition}
Given a non-empty set $Q$, denote with $\bQ^i$ the set of all $i^{\rm th}$ tuples with elements from $Q$. Introduce the concatenation of two tuples composed of elements (symbols) $\alpha_1,\ldots,\alpha_i$, and $\beta_1,\ldots,\beta_j$ from $Q$ as the mapping $\circledcirc: \bQ^i \times \bQ^j \rightarrow \bQ^{i+j}$ with the following property: 
$$
\big{(} \alpha_1 , \alpha_2 , \ldots, \alpha_i \big{)} \circledcirc \big{(} \beta_1, \beta_2, \ldots, \beta_j \big{)} = \big{(} \alpha_1, \alpha_2, \ldots \alpha_i, \beta_1, \beta_2, \ldots \beta_j \big{)}, \ \ 
$$

\end{definition}

\begin{remark}
{\rm In the following we denote the $\ell^{th}$ element of the ordered set $\bmu_{j}^{(i)}$ with $\bmu_j^{(i)}(\ell)$ (where $j \in Q$ and $i \geqslant 1$). For instance, $\bmu_1^{(2)}(3) := \big{(} \mu_1^{(2)},\mu_4^{(2)},\mu_5^{(2)} \big{)}$.
 For simplicity, use the notation $\bH_{1,2,1}(\mu_1^{(2)},\mu_4^{(2)},\mu_5^{(2)})$ instead of $\bH_{( 1,2,1 )}\big{(} ( \mu_1^{(2)},\mu_4^{(2)},\mu_5^{(2)} ) \big{)}$.}
\end{remark}

\begin{definition}
We define the {\bf nested right multi-tuples}:
\begin{equation}\label{rmtup}
\blambda_1=\left\{\blambda_1^{(1)},\blambda_1^{(2)},\ldots,\blambda_1^{({k^\dagger})}\right\}, \ \ \blambda_2=\left\{\blambda_2^{(1)},\blambda_2^{(2)},\ldots,\blambda_2^{({{k^\dagger}})}\right\}
\end{equation}
composed of sets of right $i^{\rm th}$ tuples:
\begin{equation} \label{rtup}
\hspace*{-3mm}
\begin{array}{ll}
\blambda_1^{(i)}\!=\!\left\{\!\! \begin{array}{l}
\big{(} \lambda^{(i)}_{1} \big{)},\,\\[2mm]
\big{(} \lambda^{(i)}_3,\,\lambda^{(i)}_{2} \big{)},\,\\[2mm]
\big{(} \lambda^{(i)}_5,\,\lambda^{(i)}_{4},\,\lambda^{(i)}_{1}\big{)}, \\[1mm]
\quad\quad \vdots \\[1mm]
\big{(} \lambda^{(i)}_{2m_i-3},\,\ldots,\,\lambda^{(i)}_4,\,\lambda^{(i)}_{1} \big{)},\,\\[2mm]
\big{(} \lambda^{(i)}_{2m_i-1},\,\lambda^{(i)}_{2m_i-2},\,\ldots,\,\lambda^{(i)}_3,\,\lambda^{(i)}_{2} \big{)}\\
\end{array} \right.\!\!\!\! ,\!\!
&
\blambda_2^{(i)}\!=\!\left\{\!\! \begin{array}{l}
\big{(} \lambda^{(i)}_{2} \big{)},\,\\[2mm]
\big{(} \lambda^{(i)}_4,\,\lambda^{(i)}_{1} \big{)},\,\\[2mm]
\big{(} \lambda^{(i)}_6,\,\lambda^{(i)}_{3},\,\lambda^{(i)}_{2} \big{)}, \\[1mm]
\quad\quad \vdots \\[1mm]
\big{(} \lambda^{(i)}_{2m_i-2},\,\ldots,\,\lambda^{(i)}_3,\,\lambda^{(i)}_{2} \big{)},\,\\[2mm]
\big{(} \lambda^{(i)}_{2m_i},\,\lambda^{(i)}_{2m_i-3},\,\ldots,\,\lambda^{(i)}_4,\,\lambda^{(i)}_{1} \big{)} \\
\end{array} \right.\!\!\!\! \!\!
\end{array}
\end{equation}
That is, we construct the right tuples based on the following recurrence relations ( where $\blambda_1^{(i)} (1) = \big{(} \lambda_1^{(i)} \big{)}$ and $\blambda_2^{(i)} (1) = \big{(} \lambda_2^{(i)} \big{)}$)
\begin{equation}
\blambda_1^{(i)} (g) = \big{(} \lambda_{2g-1}^{(i)} \big{)} \circledcirc \blambda_2^{(i)} (g-1), \ \ \ \blambda_2^{(i)} (g) = \big{(} \lambda_{2 g }^{(i)} \big{)}  \circledcirc \blambda_1^{(i)} (g-1),
\vspace{-2mm}
\end{equation}
for $\lambda_{2g-1}^{(i)} \in \mathbb{C}, \ 1 < g \leqslant m_i, \ i = 1,\ldots,k^{\dagger},  \ m_i \geqslant 1$ so that the equality $m_1+\,\cdots\,+m_{{k^\dagger}}={k}$ holds.\\

\vspace{-1mm}
\noindent
Also, introduce the {\bf nested left multi-tuples}
\begin{equation}\label{lmtup}
\bmu_1=\left\{\bmu_1^{(1)},\bmu_1^{(2)},\ldots,\bmu_1^{({{\ell^\dagger}})}\right\} , \ \ \bmu_2=\left\{\bmu_2^{(1)},\bmu_2^{(2)},\ldots,\bmu_2^{({{\ell^\dagger}})}\right\}
\end{equation}
composed of sets of left $j^{\rm th}$ tuples
\begin{equation} \label{ltup}
\hspace*{-3mm}
\begin{array}{ll}
 \!\bmu_1^{(j)}\!=\!\left\{\!\! \begin{array}{l}
\big{(} \mu^{(j)}_{1} \big{)},\,\\[2mm]
\big{(} \mu^{(j)}_2,\,\mu^{(j)}_{3} \big{)},\,\\[2mm]
\big{(} \mu^{(j)}_1,\,\mu^{(j)}_{4},\mu^{(j)}_{5} \big{)},\,\\[1mm]
\quad\quad \vdots \\[1mm]
\big{(} \mu^{(j)}_1,\,\mu^{(j)}_{4},\,\ldots,\,\ \mu^{(j)}_{2p_j-3} \big{)} ,\,\\[2mm]
\big{(} \mu^{(j)}_2,\,\mu^{(j)}_{3},\,\ldots,\,\mu^{(j)}_{2p_j-2},\,\mu^{(j)}_{2p_j-1} \big{)} \,\\
\end{array} \right.
&
 \!\bmu_2^{(j)}\!=\!\left\{\!\! \begin{array}{l}
\big{(} \mu^{(j)}_{2} \big{)},\,\\[2mm]
\big{(} \mu^{(j)}_1,\,\mu^{(j)}_{4} \big{)},\,\\[2mm]
\big{(} \mu^{(j)}_2,\,\mu^{(j)}_{3},\mu^{(j)}_{6} \big{)}, \,\\[1mm]
\quad\quad \vdots \\[1mm]
\big{(} \mu^{(j)}_2,\,\mu^{(j)}_{3},\,\ldots,\,\ \mu^{(j)}_{2p_j-2}\big{)},\,\\[2mm]
\big{(} \mu^{(j)}_1,\,\mu^{(j)}_{4},\,\ldots,\,\mu^{(j)}_{2p_j-3},\,\mu^{(j)}_{2p_j} \big{)} \,\\
\end{array} \right.
\end{array}\!\!\!\!\!
\end{equation}
That is, we construct the left tuples based on the following recurrence relations (where $\bmu_1^{(j)} (1) = \big{(} \mu_1^{(j)} \big{)}$ and $\bmu_2^{(j)} (1) = \big{(} \mu_2^{(j)} \big{)}$)
\begin{equation}
\bmu_1^{(j)} (h) =   \bmu_2^{(j)} (h-1) \circledcirc \big{(} \mu_{2 h-1}^{(j)} \big{)} , \ \ \ \bmu_2^{(j)} (h) =   \bmu_1^{(j)} (h-1) \circledcirc \big{(} \mu_{2 h}^{(j)} \big{)},
\vspace{-2mm}
\end{equation}
for $\mu_{2h-1}^{(j)} \in \mathbb{C}, \ 1 < h \leqslant p_j, \ j = 1,\ldots,\ell^{\dagger},  \ p_j \geqslant 1$ so that $p_1+\,\cdots\,+p_{{\ell^\dagger}}={\ell}$.
\end{definition}
\noindent
Given that the following conditions are satisfied for all $i = 1,\ldots,k^{\dagger}$ and $g = 1,\ldots,m_i$,
\begin{equation}\label{cond_lam}
\lambda_{2g-1}^{(i)} \notin {\rm eig}(\bA_1,\bE_1), \ \ \ \lambda_{2g}^{(i)} \notin {\rm eig}(\bA_2,\bE_2)
\end{equation}
associate the following matrices to the set of right tuples in (\ref{rtup}), as
$$
\cR_1^{(i)}= \left[ \bPhi_1({\lambda^{(i)}_{1}})\,\bB_1 ,~ \bPhi_1({\lambda^{(i)}_{3}})\, \bK_2 \, \bPhi_2({\lambda^{(i)}_{2}}) \, \bB_2, ~ \ldots ,~\ \bPhi_1({\lambda^{(i)}_{2m_i-1}})\,\bK_2 \cdots \bK_1\,\bPhi_1({\lambda^{(i)}_{3}})\,\bK_2\, 
\bPhi_2({\lambda^{(i)}_{2}})\,\bB_2 \right],
\vspace{-1mm}
$$
$$
\cR_2^{(i)}= \left[ \bPhi_2({\lambda^{(i)}_{2}})\,\bB_2 ,~ \bPhi_2({\lambda^{(i)}_{4}}) \, \bK_1 \, \bPhi_1({\lambda^{(i)}_{1}})\,\bB_1, ~ \ldots ,~\ \bPhi_2({\lambda^{(i)}_{2m_i}})\,\bK_1 \cdots \bK_2\,\bPhi_2({\lambda^{(i)}_{2}})\,\bK_1\, 
\bPhi_1({\lambda^{(i)}_{1}})\,\bB_1 \right],
$$
where $i = 1,\ldots,{k^\dagger}$ and $ \cR_q^{(i)} \in \IC^{n_q\times m_i} $ is attached to $\bLambda_q^{(i)}$ for $q \in \{1,2\}$. The matrices:
\begin{equation} \label{reach}
\cR_1=\left[\cR_1^{(1)},~\cR_1^{(2)},~\ldots,~\cR_1^{({k^\dagger})}\right]\in\IC^{n_1 \times {k}}, \ \  \cR_2=\left[\cR_2^{(1)},~\cR_2^{(2)},~\ldots,~\cR_2^{({k^\dagger})}\right]\in\IC^{n_2 \times {k}}.
\end{equation}
are defined as the {\bf generalized controllability matrix} 
of the LSS system $\Si$, associated with the right multi-tuple $\blambda$.
Similarly, assuming that the following conditions are satisfied for all $j = 1,\ldots,\ell^{\dagger}$ and $h = 1,\ldots,p_j$,
\begin{equation}\label{cond_mu}
\mu_{2h-1}^{(j)} \notin {\rm eig}(\bA_1,\bE_1), \ \ \  \mu_{2h}^{(j)} \notin {\rm eig}(\bA_2,\bE_2)
\end{equation}
associate the following matrices to the set of right tuples in (\ref{ltup}), as
$$
\cO_1^{(j)}= \left[ \begin{array}{c} \bC_1\,\bPhi_1({\mu^{(j)}_{1}}) \\ \bC_2\,\bPhi_2({\mu^{(j)}_{2}}) \, \bK_1 \,\bPhi_1({\mu^{(j)}_{3}})   \\ \vdots \\ \bC_2\,\bPhi_2({\mu^{(j)}_{2}})\,\bK_1\,\bPhi_1({\mu^{(j)}_{3}})\,\bK_2\, \cdots\,\bK_1\,
\bPhi_1({\mu^{(j)}_{2p_j-1}}) \end{array} \right] \in\IC^{p_j \times n_1}   , ~~j=1,\ldots,{\ell^\dagger} , 
$$
$$
\cO_2^{(j)}= \left[ \begin{array}{c}  \bC_2\,\bPhi_2({\mu^{(j)}_{1}}) \\ \bC_1\,\bPhi_1({\mu^{(j)}_{1}}) \, \bK_2 \,\bPhi_2({\mu^{(j)}_{4}})   \\  \vdots \\ \bC_1\,\bPhi_1({\mu^{(j)}_{1}})\,\bK_2\,\bPhi_2({\mu^{(j)}_{2}})\,\bK_1\, \cdots\,\bK_2\,
\bPhi_2({\mu^{(j)}_{2p_j}}) \end{array} \right] \in\IC^{p_j \times n_2}   , ~~j=1,\ldots,{\ell^\dagger} , 
$$
and the {\bf generalized observability matrices}:
\begin{equation} \label{obs}
 \cO_1=\left[\begin{array}{c}
\cO_1^{(1)}\\\vdots\\[1mm]\cO_1^{({\ell^\dagger})}\end{array}\right]\in\IC^{\ell \times n_1}, \ \  \cO_2=\left[\begin{array}{c}
\cO_2^{(1)}\\\vdots\\[1mm]\cO_2^{({\ell^\dagger})}\end{array}\right]\in\IC^{\ell \times n_2}.
\end{equation}

\vspace{-2mm}

\begin{definition}
For $\nu \in \{1,2\}$, let $\bQ^{\nu,+}$ and $\bQ^{+,\nu}$ be the ordered sets containing all tuples that can be constructed with symbols from the alphabet $Q = \{1,2\}$ and that \textbf{start} (and respectively \textbf{end}) with the symbol $\nu$. Also, no two consecutive characters are allowed to be the same. Hence, explicitly write the new introduced sets as follows:
\vspace{-2mm}
\begin{eqnarray}
&& \bQ^{1,+} = \{( 1 ), ( 1,2), (1,2,1), \ldots \}, \ \ \bQ^{2,+} = \{(2),(2,1),(2,1,2),\ldots \}, \\
&& \bQ^{+,1} = \{ (1),(2,1),(1,2,1),\ldots \}, \ \ \bQ^{+,2} = \{(2),(1,2),(2,1,2),\ldots \}
\vspace{-1mm}
\end{eqnarray}
\end{definition}

\vspace{-2mm}

\begin{remark}
{\rm In the following we denote the $\ell^{th}$ element of the ordered set $\bQ^{\nu,+}$ with $\bQ^{\nu,+}(\ell)$. For example, one writes $\bQ^{1,+}(4) := (1,2,1,2)$. Moreover, we have $\bQ^{+,2}(3) \circledcirc \bQ^{1,+}(2) = (2,1,2,1,2)$. }
\end{remark}
The compact notation $\bH_{\bQ^{+,1}}(\bmu_1(2))$ is used instead of $\bH_{2,1}(\mu_2,\mu_3)$, where $\bmu_1(2) :=  \big{(} \mu_2,\mu_3 \big{)}$


\begin{definition}
 Let the $i^{th}$ unit vector be denoted with $\bfe_i=[0~\ldots,1,\ldots,0]^T\in\IR^{ k}$. In some contexts we may use the alternative notation $\bfe_{i,k}$ to stress the fact the vector has dimension $k$. 
 
Also let $\bfz_{k,\ell} \in \mathbb{R}^{k \times \ell}$ be an all zero matrix. Hence, use the notation $\bfz_k = [0,0,\ldots,0]^T\in\IR^{ k}$ for the zero valued vector of size $k$.
 \end{definition}

In the sequel, denote with $\hat{\bH}$ the {\it generalized transfer functions} corresponding to a LSS $\hat{\Si}$.

\begin{definition}
We say that a LSS $\hat{\Si} = (k,k, \{(\hat{\bE}_i,\hat{\bA}_i, \hat{\bB}_i, \hat{\bC}_i\}_{i=1}^2,\{ \hat{\bK}_{i,j}\}_{i,j = 1}^2, \bfz)$ matches the \\[1mm] data associated with the right tuples $\{\blambda_{a}^{(1)},\ldots,\blambda_{a}^{(k^{\dagger})}\}$ as well as left tuples $\{\bmu_{b}^{(1)},\ldots,\bmu_{b}^{(\ell^{\dagger})}\}, \ a,b =$  \\[1mm] $1,2$ and corresponding to the original LSS $\Si = (n_1,n_2, \{(\bE_i,\bA_i, \bB_i, \bC_i)\}_{i=1}^2,\{ \bK_{i,j} \}_{i,j=1}^2, \bfz)$, if the following $2(k^2+2k)$ relations
\begin{equation}
\begin{cases} 
 \bH_{\bQ^{+,1}(h)}(\bmu_1^{(j)}(h)) =\hat{\bH}_{\bQ^{+,1}(h)}(\bmu_1^{(j)}(h)), \ \  \bH_{\bQ^{+,2}(h)}(\bmu_2^{(j)}(h)) =\hat{\bH}_{\bQ^{+,2}(h)}(\bmu_2^{(j)}(h)),   \\
 \bH_{\bQ^{1,+}(g)}(\blambda_1^{(i)}(g)) =\hat{\bH}_{\bQ^{1,+}(g)}(\blambda_1^{(i)}(g)), \ \  
\bH_{\bQ^{2,+}(g)}(\blambda_2^{(i)}(g)) =\hat{\bH}_{\bQ^{2,+}(g)}(\blambda_2^{(i)}(g)) ,  \\[1mm]
 \bH_{\bQ^{+,1}(h) \circledcirc \bQ^{2,+}(g)}(\bmu_1^{(j)}(h) \circledcirc \blambda_2^{(i)}(g)) =\hat{\bH}_{\bQ^{+,1}(h) \circledcirc \bQ^{2,+}(g)}(\bmu_1^{(j)}(h) \circledcirc \blambda_2^{(i)}(g)) \\[1mm]
 \bH_{\bQ^{+,2}(h) \circledcirc \bQ^{1,+}(g)}(\bmu_2^{(j)}(h) \circledcirc \blambda_1^{(i)}(g)) =\hat{\bH}_{\bQ^{+,2}(h) \circledcirc \bQ^{1,+}(g)}(\bmu_2^{(j)}(h) \circledcirc \blambda_1^{(i)}(g)) 
\end{cases}
\end{equation}
hold for $j = 1,\ldots,k^{\dagger}, \ h = 1,\ldots,p_j$ and $i = 1,\ldots,k^{\dagger}, \ g = 1,\ldots,m_i$, where
$$
p_1+p_2+\ldots p_{k^{\dagger}} = m_1+m_2+\ldots m_{k^{\dagger}} = k.
$$
\end{definition}

\noindent
The following lemma extends the rational interpolation idea for linear systems approximation to the linear switched system case. 

\vspace{-2mm}

\begin{lemma} {\bf Interpolation of LSS}. 
Let $\Si = (n_1,n_2,\{(\bE_i,\bA_i, \bB_i, \bC_i)\}_{i=1}^2,\{ \bK_{i,j}  \}_{i,j =1}^2, \bfz)$ be a LSS of order $(n_1,n_2)$. An  order $k$  reduced LSS $\hat{\Si} = (k,k,\{(\hat{\bE}_i,\hat{\bA}_i, \hat{\bB}_i, \hat{\bC}_i)\}_{i=1}^2,\{ \hat{\bK}_{i,j}  \}_{i,j =1}^2, \bfz)$ is constructed using the projection matrices chosen as in (\ref{reach}) and (\ref{obs}) for $\ell=k$, i.e. 
$$
\bX_1=\cR_1, \ \ \bX_2=\cR_2 \ \text{and} \ \bY_1^T=\cO_1, \ \ \bY_2^T=\cO_2 
\vspace{-1mm}
$$
 Additionally assume ${\rm rank}(\cR_i) = {\rm rank}(\cO_i) = k, \ i \in \{1,2\}$. The reduced matrices corresponding to the $I^{st}$ subsystem $\hat{\Si}_1$ are computed as,
\begin{equation} \label{eq:interp1sys}
\hat{\bE}_1=\bY_1^T\bE_1\bX_1, ~~\hat{\bA}_1=\bY_1^T\bA_1\bX_1,~~ \hat{\bB}_1=\bY_1^T\bB_1,~~ \hat{\bC}_1=\bC_1\bX_1,~~ \hat{\bK}_1=\bY_2^T\bK_1\bX_1, 
\end{equation}
\vspace{-1mm}
while the reduced matrices corresponding to the $II^{\rm nd}$ subsystem $\hat{\Si}_2$ can also be computed as,
\begin{equation} \label{eq:interp2sys}
\hat{\bE}_2=\bY_2^T\bE_2\bX_2, ~~\hat{\bA}_2=\bY_2^T\bA_2\bX_2,~~ \hat{\bB}_2=\bY_2^T\bB_2,~~ \hat{\bC}_2=\bC_2\bX_2,~~ \hat{\bK}_2=\bY_1^T\bK_2\bX_2. 
\end{equation}
\noindent
It follows that the reduced-order system $\hat{\Si}$ matches the
data of the system $\Si$ (as it was previously introduced in Definition 4.6). \end{lemma}

\textbf{Proof of Lemma 4.1} 
For simplicity, assume that we have one set of right multi-tuples, and one set of left multi-tuples with the same number of interpolation points $k$ for each mode. Moreover let $k$ be an even positive number. For the first mode, write down the interpolation nodes as follows:
\begin{equation} \label{eq:int_pt_mode1}
\begin{cases} \begin{array}{l}
\blambda_1=\left\{\big{(} \lambda_1 \big{)},\big{(} \lambda_3,\lambda_2 \big{)} ,~\ldots,~ \big{(} \lambda_{2k-1},\cdots,\lambda_3,\lambda_2 \big{)}
\right\},\\[1mm]
\bmu_1=\left\{\big{(} \mu_1 \big{)},\big{(} \mu_2,\mu_3 \big{)},~\ldots,~\big{(} \mu_2,\mu_3,\cdots,\mu_{2k-1} \big{)}
\right\}.
\end{array} \end{cases}
\end{equation}
For the second mode, write down the interpolation nodes as follows:
\begin{equation} \label{eq:int_pt_mode2}
\begin{cases}
\begin{array}{l}
\blambda_2=\left\{\big{(} \lambda_2 \big{)}, \big{(} \lambda_4,\lambda_1 \big{)},~\ldots,~ \big{(} \lambda_{2k},\cdots,\lambda_2,\lambda_1 \big{)}
\right\},\\[1mm]
\bmu_2=\left\{ \big{(} \mu_2 \big{)}, \big{(} \mu_1,\mu_4 \big{)},~\ldots,~\big{(} \mu_1,\mu_4,\cdots,\mu_{2k} \big{)}
\right\}.
\end{array} \end{cases}
\end{equation}
This corresponds to the case for which $l = k, l^\dagger = k^\dagger = 1$ and $m_1 = p_1 = k$. It follows that the interpolation conditions stated in Definition 4.6, can be rewritten by taking into account the aforementioned simplification as, 
\begin{eqnarray} 
& \underline{\mathbf{2k} \ \textbf{conditions:}} \ \begin{cases} \bH_{\bQ^{+,1}(j)}(\bmu_1(j)) =\hat{\bH}_{\bQ^{+,1}(j)}(\bmu_1(j)) \\  \bH_{\bQ^{+,2}(j)}(\bmu_2(j)) =\hat{\bH}_{\bQ^{+,2}(j)}(\bmu_2(j))
 \end{cases}, \ \ j \in \{1,\ldots,k\} \label{eq:interp_cond1} \\
& \underline{\mathbf{2k} \ \textbf{conditions:}} \ \begin{cases} \bH_{\bQ^{1,+}(i)}(\blambda_1(i)) =\hat{\bH}_{\bQ^{1,+}(i)}(\blambda_1(i)) \label{eq:interp_cond2} \\ 
\bH_{\bQ^{2,+}(i)}(\blambda_2(i)) =\hat{\bH}_{\bQ^{2,+}(i)}(\blambda_2(i)) \end{cases} , \ \ i \in \{1,\ldots,k\}  \\[1mm]
& \underline{\mathbf{k^2} \ \textbf{conditions:}} \begin{cases} \bH_{\bQ^{+,1}(j) \circledcirc \bQ^{2,+}(i)}(\bmu_1(j) \circledcirc \blambda_2(i)) =\hat{\bH}_{\bQ^{+,1}(j) \circledcirc \bQ^{2,+}(i)}(\bmu_1(j) \circledcirc \blambda_2(i)), \end{cases} \label{eq:interp_cond3} \\
& \underline{\mathbf{k^2} \ \textbf{conditions:}} \begin{cases} \bH_{\bQ^{+,2}(j) \circledcirc \bQ^{1,+}(i)}(\bmu_2(j)\circledcirc \blambda_1(i)) =\hat{\bH}_{\bQ^{+,2}(j) \circledcirc \bQ^{1,+}(i)}(\bmu_2(j) \circledcirc \blambda_1(i))  \end{cases}& \label{eq:interp_cond4}
\end{eqnarray}

With the assumptions in (\ref{eq:int_pt_mode1}) and (\ref{eq:int_pt_mode2}), it follows that the associated generalized controllability and observability matrices defined previously in (\ref{reach}) and (\ref{obs}), are rewritten as:
$$
\begin{array}{l}
\cR_1=\left[~\bPhi_1(\lambda_1)\bB_1,~\bPhi_1(\lambda_3)\bK_2 \bPhi_2(\lambda_2)\bB_2,~\ldots,~
\bPhi_1(\lambda_{2k-1})\bK_2 \, \cdots \,\bK_2 \bPhi_2(\lambda_2)\bB_2
\right]\in\IC^{n\times k},\\[2mm]
\cR_2=\left[~\bPhi_2(\lambda_2)\bB_2,~\bPhi_2(\lambda_4)\bK_1 \bPhi_1(\lambda_1)\bB_1,~\ldots,~
\bPhi_2(\lambda_{2k})\bK_1 \, \cdots \,\bK_1 \bPhi_1(\lambda_1)\bB_1
\right]\in\IC^{n\times k}, \\ [2mm]
\cO_1=\left[\begin{array}{l}
\bC_1\bPhi_1(\mu_1)\\
\bC_2\bPhi_2(\mu_2)\bK_1\bPhi_1(\mu_3)\\
\qquad\vdots\\
\bC_2\bPhi_2(\mu_2)\bK_1 \bPhi_1(\mu_3) \cdots \bK_1\bPhi_1(\mu_{2k-1})\\
\end{array}\right],  \ \cO_2=\left[\begin{array}{l}
\bC_2\bPhi_2(\mu_2)\\
\bC_1\bPhi_1(\mu_1)\bK_2\bPhi_2(\mu_4)\\
\qquad\vdots\\
\bC_1\bPhi_1(\mu_1)\bK_2 \bPhi_2(\mu_4) \cdots  \bK_2\bPhi_2(\mu_{2k})\\
\end{array}\right].
\end{array}
$$
with both $\cO_1, \ \cO_2  \in\IC^{k\times n}$. Additionally, introduce the notation $\hat{\bPhi}_i(s) = (s \hat{\bE}-\hat{\bA})^{-1}$.\\
From (\ref{eq:interp1sys}) and (\ref{eq:interp2sys}), using that $\bX_i = \cR_i$ for $i = 1,2$, it readily follows that: 
$$
{\bf (a)}~~\hat\bPhi_1(\lambda_1)\,\hat\bB_1=\bfe_1\quad\mbox{and}\quad
{\bf (b)}~~ \hat\bPhi_1 (\lambda_{2i-1})\,\hat\bK_2\,\bfe_{i-1}=\,\bfe_i,~i=2,\ldots,k,
$$
$$
{\bf (c)}~~\hat\bPhi_2(\lambda_2)\,\hat\bB_2=\bfe_1\quad\mbox{and}\quad
{\bf (d)}~~ \hat\bPhi_2 (\lambda_{2i})\,\hat\bK_1\,\bfe_{i-1}=\,\bfe_i,~i=2,\ldots,k .
$$
These equalities imply the right-hand conditions in (\ref{eq:interp_cond2}).
Similarly, from (\ref{eq:interp1sys}) and (\ref{eq:interp2sys}), using that $\bY_j^T = \cO_j$ for $j = 1,2$, it follows that:
$$
{\bf (e)}~~\bC_1\,\hat\bPhi_1(\mu_1)=\bfe_1^T\quad\mbox{and}\quad
{\bf (f)}~~ \bfe_{j-1}^T\bK_2\hat\bPhi_2 (\mu_{2j})\,=\,\bfe_j^T,~j=2,\ldots,k ,
$$
$$
{\bf (g)}~~\bC_2\,\hat\bPhi_2(\mu_2)=\bfe_1^T\quad\mbox{and}\quad
{\bf (h)}~~ \bfe_{j-1}^T\bK_1\hat\bPhi_1 (\mu_{2j-1})\,=\,\bfe_j^T,~j=2,\ldots,k ,
$$
which imply left-hand conditions in (\ref{eq:interp_cond1}). Finally, with $\bX=\cR$, $\bY^T=\cO$, and
combining {\bf (a)-(h)}, all interpolation conditions in (\ref{eq:interp_cond3}) and (\ref{eq:interp_cond4}) are hence satisfied. \sq

\vspace{-3mm}

\begin{remark}
{\rm For instance, in Example 4.2, the conditions stated in (\ref{eq_ex_cond}) are satisfied.}
\end{remark}

\vspace{-6mm}

\subsubsection{Sylvester equations for $\cO$ and $\cR$}

The generalized controllabilty and observability matrices satisfy Sylvester equations. To state the corresponding result we need to define the following quantities. First introduce the matrices
\begin{equation} \label{RMat}
\bR=\left[\bfe^T_{1,m_1}~\cdots~\bfe^T_{1,m_{k^\dagger}}\right] \in\IR^{1 \times k},~~
\bL^T=\left[\bfe^T_{1,p_1}~\cdots~\bfe^T_{1,p_{\ell^\dagger}}\right] \in\IR^{1 \times \ell},
\end{equation}
and the block-shift matrices
\footnotesize\small
\begin{equation} \label{SMat}
\hspace*{-4mm}
\left\{\begin{array}{ll}
\bS_\bR&=~\mbox{blkdiag}\,\left[\bJ_{m_1},~\ldots,~\bJ_{m_{k^\dagger}}\right],\\[1mm]
\bS_\bL&=~\mbox{blkdiag}\,\left[\bJ_{p_1}^T,~\ldots,~\bJ_{p_{\ell^\dagger}}^T\right] .\\ 
\end{array}\right.  \ \ \text{where} \ \ 
\bJ_{u} = \left(\begin{array}{cccc} 
0 & 1 & \cdots   & 0\\  
\vdots & \vdots & \ddots & \vdots  \\  
0 & 0 & \cdots   & 1 \\ 
0 & 0 & \cdots   & 0  \end{array}\right)
\in\IR^{u \times u}~
\end{equation} 
\normalsize
Finally we arrange the left interpolation points in the diagonal matrices as,
\begin{equation} \label{MLam}
\bM_1 = \mbox{blkdiag}\,[\bM_1^{(1)},~~\bM_1^{(2)},~\ldots,~\bM_1^{(\ell^\dagger)}], \ \ 
 \bM_2 = \mbox{blkdiag}\,[\bM_2^{(1)},~~\bM_2^{(2)},~\ldots,~\bM_2^{(\ell^\dagger)}],
\end{equation} 
where $\bM_1^{(j)} = \mbox{diag}\,[\mu_1^{(j)},~~\mu_3^{(j)},~\ldots,~\mu_{2p_j-1}^{(j)}]$ \ and \ $\bM_2^{(j)} = \mbox{diag}\,[\mu_2^{(j)},~~\mu_4^{(j)},~\ldots,~\mu_{2p_j}^{(j)}]$; we used the MATLAB notation 'blkdiag' which outputs a block diagonal matrix with each input entry as a block.
Also arrange the right interpolation points in the diagonal matrices:
\begin{equation} \label{LLam}
\bLambda_1 = \mbox{blkdiag}\,[\bLambda_1^{(1)},~~  \bLambda_1^{(2)},~\ldots,~ \bLambda_1^{(\ell^\dagger)}], \ \ \bLambda_2 = \mbox{blkdiag}\,[\bLambda_2^{(1)},~~  \bLambda_2^{(2)},~\ldots,~ \bLambda_2^{(\ell^\dagger)}],
\end{equation}
where $\bLambda_1^{(i)} = \mbox{diag}\,[\lambda_1^{(i)},~~\lambda_3^{(i)},~\ldots,~\lambda_{2m_i-1}^{(i)}]$ \ and \ $\bLambda_2^{(i)} = \mbox{diag}\,[\lambda_2^{(i)},~~\lambda_4^{(i)},~\ldots,~\lambda_{2m_i}^{(i)}]$.

In the definitions above, i.e. (\ref{RMat})-(\ref{LLam}) we analyzed the general case, i.e., the assumptions made in (\ref{eq:int_pt_mode1})-(\ref{eq:int_pt_mode2}) were no longer valid. We are now ready to state the following result.
\begin{lemma} \label{lemma_sylv1}
The generalized controllability matrices $\cR_1, \cR_2$ defined by (\ref{reach}), satisfy the following Sylvester equations:
\begin{equation}\label{eq:reach} 
\begin{cases}  \bA_1 \mathcal{R}_1  +  \bK_2 \mathcal{R}_2 \bS_{\bR} + \bB_1 \bR = \bE_1 \mathcal{R}_1 \bLambda_1,  \\  \bA_2 \mathcal{R}_2 +  \bK_1 \mathcal{R}_1 \bS_{\bR} +  \bB_2 \bR = \bE_2 \mathcal{R}_2 \bLambda_2.  \end{cases}
\end{equation}
\end{lemma}

\textbf{Proof of Lemma 4.2}
 Assume again, for simplicity of the proof, that the assumptions made in (\ref{eq:int_pt_mode1})-(\ref{eq:int_pt_mode2}) are valid. Hence, we have one set of right multi-tuples for each of the two modes with same number of interpolation points $k$ (with $k$ even).
 Multiplying the first equation in (\ref{eq:reach}) on the right with the first unit vector $\bfe_1$  we obtain:
\begin{equation}\label{init_cond1}
\bA_1 \cR_1^{(1)}+ \bB_1  = \lambda_1 \bE_1 \cR_1^{(1)} \Rightarrow \cR_1^{(1)} = (\lambda_1 \bE_1- \bA_1)^{-1}\bB_1 = \bPhi_1(\lambda_1) \bB_1.
\end{equation}
where $\cR_i^{(j)}$ is the $j^{th}$ column of $\cR_i$ (with $j \leqslant k$ and $i \in \{1,2\}$). 
Thus the first column of the matrix which is the solution of the first equation in (\ref{eq:reach}) is indeed 
equal to the first column of the generalized controllability matrix $\cR_1$. Multiplying the second equation in (\ref{eq:reach}) on the right with the first unit vector $\bfe_1$  we obtain:
\begin{equation}\label{init_cond2}
\bA_2 \cR_2^{(1)}+ \bB_2  = \lambda_2 \bE_2 \cR_2^{(1)} \Rightarrow \cR_2^{(1)} = (\lambda_2 \bE_2- \bA_2)^{-1}\bB_2 = \bPhi_2(\lambda_2) \bB_2.
\end{equation}
Thus the first column of the matrix which is the solution of the second equation in (\ref{eq:reach}) is indeed 
equal to the first column of the generalized controllability matrix $\cR_2$. Multiplying first equation in (\ref{eq:reach}) on the right with the $j^{th}$ unit vector $\bfe_j$, we obtain:
\begin{equation}\label{reach01}
\bA_1 \cR_1^{(j)}+ \bK_2 \cR_{2}^{(j-1)}= \lambda_{2j-1} \bE_1 \cR_1^{(j)} ~\Rightarrow~
\cR_1^{(j)} = (\lambda_{2j-1} \bE_1- \bA_1)^{-1}\bK_2 \cR_2^{(j-1)}
\end{equation}
Multiplying second equation in (\ref{eq:reach}) on the right with the $j^{th}$ unit vector $\bfe_j$, we obtain:
\begin{equation}\label{reach02}
\bA_2 \cR_2^{(j)}+ \bK_1 \cR_{1}^{(j-1)}= \lambda_{2j} \bE_2 \cR_2^{(j)} ~\Rightarrow~
\cR_2^{(j)} = (\lambda_{2j} \bE_2- \bA_2)^{-1}\bK_1 \cR_1^{(j-1)} 
\end{equation}
From (\ref{reach01}) and (\ref{reach02}) we write the following linear recursive system of equations:
\begin{equation}
\begin{cases} \cR_1^{(j)} = \bPhi_1(\lambda_{2j-1}) \bK_2 \cR_2^{(j-1)} \\ \cR_2^{(j)} = \bPhi_2(\lambda_{2j}) \bK_1 \cR_1^{(j-1)} \end{cases} 
\end{equation}
with initial conditions (\ref{init_cond1}) and (\ref{init_cond2}). Hence, by solving the coupled system of equations, we indeed conclude that $\cR_1$ and $\cR_2$ matrices satisfying (\ref{eq:reach}) are the generalized controllability matrices defined in (\ref{reach}) (for this particular choice of the right interpolation points).
This proof can be nevertheless adapted from the simplified case  in  (\ref{eq:int_pt_mode1})-(\ref{eq:int_pt_mode2}) to the more general case treated in Definition 4.3. \sq

\begin{remark}
{\rm The generalized Sylvester equations in (\ref{eq:reach}) can be compactly written as only one generalized Sylvester equation in the following way
\begin{equation}
\bA^\bD \cR^\bD+ \bK^{\scriptsize\reflectbox{\bD}} \cR^\bD \bS_{\bR}^{\scriptsize\reflectbox{\bD}}+ \bB^\bD \bR^\bD = \bE^\bD \cR^\bD \bLambda^\bD
\end{equation}
where $\bX^\bD = \left( \begin{array}{cc}
\bX_1 & \bfz \\ \bfz & \bX_2
\end{array} \right)$ for $\bX \in \{\cR,\bA,\bB,\bE,\bR,\bLambda\}$ and  $\bK^{\scriptsize\reflectbox{\bD}} = \left( \begin{array}{cc}
 \bfz & \bK_2 \\  \bK_1 & \bfz
\end{array} \right)$, ${\bS}_{\bR}^{\scriptsize\reflectbox{\bD}} = \left( \begin{array}{cc}
 \bfz & \bS_{\bR}  \\  \bS_{\bR}  & \bfz
\end{array} \right)$.
}
\end{remark}

\begin{proposition}
The pair of generalized Sylvester equations in (\ref{eq:reach}) has unique solutions if the right interpolation points are chosen so that the Sylvester operator
\vspace{-1mm}
$$
\cL_\cR=\bI\otimes \bA^{\bD} - \bLambda ^{\bD} \!\otimes \bE^{\bD} + \Big{(} \bS_\bR^{\scriptsize\reflectbox{\bD}} \Big{)}^T \!\otimes \bK^{\scriptsize\reflectbox{\bD}}  ,~~~
\vspace{-1mm}
$$
is invertible, i.e. have no zero eigenvalues (where $\otimes$ denotes the Kronecker product).
\end{proposition}

\begin{remark}
The motivation behind this subsection is closely related to building parametrized reduced order models. The idea is that, one can also use only one sided interpolation conditions (either left as in \ref{eq:interp_cond1} or right as in \ref{eq:interp_cond2}) to reduce the original LSS. The other degrees of freedom are given by free parameters.
Further development of these Sylvester equation was studied in \cite{agi16} (in Section 4.4) for the case of bilinear systems. 
\end{remark}

\begin{lemma} \label{lemma_sylv2}
The generalized  observability matrices $\cO_1$ and $\cO_2$ defined by (\ref{obs}) satisfy the following generalized Sylvester equations:
\begin{equation} \label{eq:obs}
\begin{cases}  \mathcal{O}_1 \bA_1  +  \bS_{\bL} \mathcal{O}_2 \bK_1 +  \bL \bC_1 = \bM_1 \mathcal{O}_1 \bE_1  \\   \mathcal{O}_2 \bA_2  +  \bS_{\bL} \mathcal{O}_1  \bK_2 + \bL \bC_2  = \bM_2 \mathcal{O}_2 \bE_2  \end{cases}
\end{equation}
\vspace{-5mm}
\end{lemma}
\textbf{Proof of Lemma 4.3}
Similar to the proof of Lemma 4.2. \sq

\begin{remark}
{\rm The generalized Sylvester equations in (\ref{eq:obs}) can be compactly written as only one generalized Sylvester equation in the following way
\begin{equation}
\cO^\bD \bA^\bD + \bS_{\bL}^{\scriptsize\reflectbox{\bD}}  \cO^\bD \bK^{\scriptsize\reflectbox{\bD}} +  \bL^\bD \bC^\bD =  \bM^\bD  \cO^\bD \bE^\bD
\end{equation}
where $\bX^\bD = \left( \begin{array}{cc}
\bX_1 & \bfz \\ \bfz & \bX_2
\end{array} \right)$ for $\bX \in \{\cO, \cC,\bL,\bM, \}$ and ${\bS}_{\bL}^{\scriptsize\reflectbox{\bD}} = \left( \begin{array}{cc}
 \bfz & \bS_{\bL}  \\  \bS_{\bL}  & \bfz
\end{array} \right)$.
}
\end{remark}

\begin{proposition}
The pair of generalized Sylvester equations in (\ref{eq:reach}) has unique solutions if the right interpolation points are chosen so that the Sylvester operator
$$
\cL_\cR= \Big{(} \bA^{\bD} \Big{)}^T \! \otimes \bI - \Big{(} \bE^{\bD} \Big{)}^T \! \otimes \bM ^{\bD}    + \Big{(} \bK^{\scriptsize\reflectbox{\bD}} \Big{)}^T \!\otimes \bS_\bL^{\scriptsize\reflectbox{\bD}}
$$
is invertible, i.e. have no zero eigenvalues (where $\otimes$ denotes the Kronecker product).
\end{proposition}

\subsection{The generalized Loewner pencil}

\begin{definition} Given a linear switched system $\Si$ as defined in (\ref{LSS_def}), let $\{\cR_1,\cR_2\}$ and $\{\cO_1,\cO_2\}$ be the controllability and observability matrices defined in (\ref{reach}), (\ref{obs}) respectively, and associated with the multi-tuples in (\ref{rmtup}), (\ref{lmtup}) respectively. The {\bf Loewner matrices} $\IL_1$ and $\IL_2$ are defined as
\vspace{-1mm}
\begin{equation} \label{loewner}
\IL_1= - \cO_1\,\bE_1\,\cR_1,~~\IL_2= - \cO_2\,\bE_2\,\cR_2 \, .
\vspace{-1mm}
\end{equation}
Additionally, the {\bf shifted Loewner} matrices $\sIL_1$ and $\sIL_2$ are defined as
\vspace{-1mm}
\begin{equation} \label{loewner_shift}
\sIL_1= - \cO_1\,\bA_1\,\cR_1,~~\sIL_2= - \cO_2\,\bA_2\,\cR_2.
\vspace{-1mm}
\end{equation}
Also define the quantities
\vspace{-1mm}
\begin{equation}\label{loewner2}
\begin{cases}\bW_1=\bC_1\,\cR_1 \\ \bW_2=\bC_2\,\cR_2 \end{cases},~~\begin{cases} \bV_1= \cO_1\,\bB_1 \\ \bV_2= \cO_2\,\bB_2 \end{cases}~~\mbox{and}~~ \begin{cases} \bXi_1=\cO_2\,\bK_1\,\cR_1 \\ \bXi_2=\cO_1\,\bK_2\,\cR_2 \end{cases} .
\vspace{-1mm}
\end{equation}
\end{definition}

\vspace{-4mm}

\begin{remark}
{\rm
In general, the Loewner matrices defined above need not have only real entries. For instance, it may happen that the samples points are purely imaginary values (on the $j \omega$ axis). In this case, we refer the readers to Section 4.3.1 in \cite{agi16}. We propose a similar method to enforce all system matrices have only real entries. In short, the sampling points have to be chosen as complex conjugate pairs; after the data is arranged into matrix format, use projection matrices as in equation (4.26) in \cite{agi16} to multiply the matrices in (\ref{loewner}), (\ref{loewner_shift}) and (\ref{loewner2}) to the left and to the right. In this way, the LSS does not change as pointed out in Definition 4.1. }
\end{remark}

\begin{remark}
{\rm Note that $\IL_k$ and $\sIL_k$ (where $k \in \{1,2\}$), as defined above, are indeed Loewner matrices, that is, they can be expressed as divided differences of appropriate transfer function values of the underlying LSS (see the following example).}
\end{remark}

\begin{example} \label{ex1} {\rm Given the LSS described by
$(\bC_j,\bE_j,\bA_j,\bB_j)$ (D = 2 and $j \in \{1,2\}$), consider the {ordered tuples} of left interpolation points:
$\big{\{} 
(\mu_{1}), \ 
(\mu_{2},\mu_{3})
\big{\}} , \ \ \big{\{} (\mu_{2}), \
(\mu_{1},\mu_{4}) \big{\}} $ and right interpolation points
$\big{\{}  (\lambda_{1}), \ (\lambda_{3},\lambda_{2}) \big{\}}, \ \big{\{} (\lambda_{2}), \ (\lambda_{4},\lambda_{1}) \big{\}}$. The associated {\it generalized observability and controllability} matrices are computed as follows 
$$
\begin{array}{l}
\cO_1=\left[ \begin{array}{c} \bC_1 \bPhi_1(\mu_1) \\  \bC_2 \bPhi_2(\mu_2) \bK_1 \bPhi_1(\mu_3)  \end{array} \right], \ \ \cO_2=\left[ \begin{array}{c} \bC_2 \bPhi_2(\mu_2) \\  \bC_1 \bPhi_3(\mu_1) \bK_2 \bPhi_2(\mu_4)  \end{array} \right] \\[5mm]
$$
$$
\mathcal{R}_1 = \left[ \begin{array}{cc} \bPhi_1(\lambda_1)\bB_1 &    \bPhi_1(\lambda_3) \bK_2 \bPhi_2(\lambda_2)\bB_2 \end{array} \right], \ \ \mathcal{R}_2 = \left[ \begin{array}{cc} \bPhi_2(\lambda_2)\bB_2 &    \bPhi_2(\lambda_4) \bK_1 \bPhi_1(\lambda_1)\bB_1 \end{array} \right]
\end{array}
$$
The projected Loewner matrices can be written in terms of the samples in the following way:
\large
\begin{eqnarray*}
&& \IL_1 = \left[ \begin{array}{cc} \frac{\bH_1(\mu_1)-\bH_1(\lambda_1)}{\mu_1-\lambda_1} & \frac{\bH_{1,2}(\mu_1,\lambda_2)-\bH_{1,2}(\lambda_3,\lambda_2)}{\mu_1-\lambda_3} \\[3mm] \frac{\bH_{2,1}(\mu_2,\mu_3)-\bH_{2,1}(\mu_2,\lambda_1)}{\mu_3-\lambda_1} & \frac{\bH_{2,1,2}(\mu_2,\mu_3,\lambda_2)-\bH_{2,1,2}(\mu_2,\lambda_3,\lambda_2)}{\mu_3-\lambda_3}
\end{array} \right] = -\mathcal{O}_1 \bE_1 \mathcal{R}_1 ,\\[1mm]
&& \IL_2 = \left[ \begin{array}{cc} \frac{\bH_2(\mu_2)-\bH_2(\lambda_2)}{\mu_2-\lambda_2} & \frac{\bH_{2,1}(\mu_2,\lambda_1)-\bH_{2,1}(\lambda_4,\lambda_1)}{\mu_2-\lambda_4} \\[3mm] \frac{\bH_{1,2}(\mu_1,\mu_4)-\bH_{1,2}(\mu_1,\lambda_2)}{\mu_4-\lambda_2} & \frac{\bH_{1,2,1}(\mu_1,\mu_4,\lambda_4)-\bH_{1,2,1}(\mu_1,\lambda_4,\lambda_1)}{\mu_4-\lambda_4}
\end{array} \right] = -\mathcal{O}_2 \bE_2 \mathcal{R}_2 
\end{eqnarray*}
\normalsize
The projected shifted Loewner matrices can also be written in terms of the samples as:
\large
\begin{eqnarray*}
&& \sIL_1 = \left[ \begin{array}{cc} \frac{\mu_1\bH_1(\mu_1)-\lambda_1\bH_1(\lambda_1)}{\mu_1-\lambda_1} & \frac{\mu_1 \bH_{1,2}(\mu_1,\lambda_2)- \lambda_3 \bH_{1,2}(\lambda_3,\lambda_2)}{\mu_1-\lambda_3} \\[3mm]  \frac{\mu_3 \bH_{2,1}(\mu_2,\mu_3)- \lambda_1 \bH_{2,1}(\mu_2,\lambda_1)}{\mu_3-\lambda_1} & \frac{\mu_3 \bH_{2,1,2}(\mu_2,\mu_3,\lambda_2)-\lambda_3 \bH_{2,1,2}(\mu_2,\lambda_3,\lambda_2)}{\mu_3-\lambda_3}
\end{array} \right] = -\mathcal{O}_1 \bA_1 \mathcal{R}_1 ,\\[1mm]
&& \sIL_2 = \left[ \begin{array}{cc} \frac{\mu_2 \bH_2(\mu_2)- \lambda_2 \bH_2(\lambda_2)}{\mu_2-\lambda_2} & \frac{\mu_2 \bH_{2,1}(\mu_2,\lambda_1)-\lambda_4 \bH_{2,1}(\lambda_4,\lambda_1)}{\mu_2-\lambda_4} \\[3mm]  \frac{\mu_4\bH_{1,2}(\mu_1,\mu_4)-\lambda_2 \bH_{1,2}(\mu_1,\lambda_2)}{\mu_4-\lambda_2} & \frac{\mu_4 \bH_{1,2,1}(\mu_1,\mu_4,\lambda_4)-\lambda_4 \bH_{1,2,1}(\mu_1,\lambda_4,\lambda_1)}{\mu_4-\lambda_4}
\end{array} \right] = -\mathcal{O}_2 \bA_2 \mathcal{R}_2 ,\\[1mm]
\end{eqnarray*}
\normalsize
The same property applies for the $\bV_i$ and $\bW_j$ vectors and $\bXi_j$ matrices:
\begin{eqnarray*}
&& \bV_1  = \left[ \begin{array}{c}
\bH_1(\mu_1) \\ \bH_{2,1}(\mu_2,\mu_3)
\end{array} \right] =  \mathcal{O}_1 \bB_1}, \ \ {\color{black}  \bV_2  = \left[ \begin{array}{c}
\bH_2(\mu_2) \\ \bH_{1,2}(\mu_1,\mu_4)
\end{array} \right] =   \mathcal{O}_2 \bB_2  ,\\[1mm]
&& \bW_1 = \left[ \begin{array}{cc}
\bH_1(\lambda_1) & \bH_{1,2}(\lambda_3,\lambda_2)
\end{array} \right] = \bC_1 \mathcal{R}_1}, \ \ {\color{black}  \bW_2 = \left[ \begin{array}{cc}
\bH_2(\lambda_2) & \bH_{2,1}(\lambda_4,\lambda_1)
\end{array} \right] = \bC_2 \mathcal{R}_2, \\[1mm]
&& \bXi_1 = \left[ \begin{array}{cc}
\bH_{2,1}(\mu_2,\lambda_1) & \bH_{2,1,2}(\mu_2,\lambda_3,\lambda_2)  \\
\bH_{1,2,1}(\mu_1,\mu_4,\lambda_1)  & \bH_{1,2,1,2}(\mu_1,\mu_4,\lambda_3,\lambda_2)  
\end{array} \right] = \mathcal{O}_2 \bK_1 \mathcal{R}_1,\\[1mm]
&& \bXi_2 = \left[ \begin{array}{cc}
\bH_{1,2}(\mu_1,\lambda_2) & \bH_{1,2,1}(\mu_1,\lambda_4,\lambda_1)  \\
\bH_{2,1,2}(\mu_2,\mu_3,\lambda_2)  & \bH_{2,1,2,1}(\mu_2,\mu_3,\lambda_4,\lambda_1)  
\end{array} \right] = \mathcal{O}_1 \bK_2 \mathcal{R}_2
\end{eqnarray*}

\normalsize
It readily follows that, given the original system $\Si$,
a reduced LSS of order two can be obtained without 
computation (matrix factorizations or solves) as:
$$
\hat\bE_k =\cO\bE\cR,~~\hat\bA=\cO\bA\cR,~~\hat\bN=\cO\bN\cR,~~\hat\bB=\cO\bB,~~\hat\bC=\bC\cR.
$$
This reduced system matches sixteen moments of the original system, namely: 
\begin{equation} \label{eq_ex_cond}
\begin{array}{rl}
\text{four of} \ ~\bH_1 / \bH_2:& \bH_1(\mu_1),~ \bH_2(\mu_2),~ \bH_1(\lambda_1), ~ \bH_2(\lambda_2), \\
\text{three of} \ ~\bH_{1,2}:& \bH_{1,2}(\mu_1,\mu_4),~\bH_{1,2}(\mu_1,\lambda_2), ~\bH_{1,2}(\lambda_3,\lambda_2),  \\
\text{three of} \ ~\bH_{2,1}:& \bH_{2,1}(\mu_2,\mu_3),~\bH_{2,1}(\mu_2,\lambda_1), ~\bH_{2,1}(\lambda_4,\lambda_1),  \\
\vdots \\
\text{one of} \ ~\bH_{1,2,1,2}:& \bH_{1,2,1,2}(\mu_1,\mu_4,\lambda_3,\lambda_2) \\
\text{one of} \ ~\bH_{2,1,2,1}:& \bH_{2,1,2,1}(\mu_2,\mu_3,\lambda_4,\lambda_1)
\end{array}
\end{equation}
i.e. in total ~$2(2k+k^2)=16$ moments are matched using this procedure.  \sq}\\
\end{example}

\subsubsection{Properties of the Loewner pencil}
We will now show that the quantities defined earlier satisfy various equations which generalize the ones in the linear or bilinear case. 

The equations that are be presented in this section are used to automatically find the Loewner and shifted Loewner matrices by means of solving Sylvester equations (instead of building the divided difference matrices from the computed samples at the sampling points).

\begin{proposition}
The Loewner matrix \,$\IL_1$ and the shifted Loewner matrix $\sIL_1$ (corresponding to mode 1) satisfy the following relations (where $\bL,\bR,\bLambda_k,\bM_k,\bS_\bL,\bS_\bR$ are given in (\ref{RMat}),(\ref{SMat}) and(\ref{MLam})): 
\begin{eqnarray} 
&& \sIL_1 = \IL_1 \bLambda_1 + \bV_1 \bR + \bXi_2 \bS_{\bR} \label{eq:LLs11} \\
&& \sIL_1 = \bM_1 \IL_1  +  \bL \bW_1  +  \bS_{\bL} \bXi_1 \label{eq:LLs12}
\end{eqnarray}
The Loewner matrix \,$\IL_2$ and the shifted Loewner matrix $\sIL_2$ (corresponding to mode 2) satisfy the following relations: 
\begin{eqnarray} 
&& \sIL_2 = \IL_2 \bLambda_2 + \bV_2 \bR +  \bXi_1 \bS_{\bR} \label{eq:LLs21} \\
&& \sIL_2 = \bM_2 \IL_2 + \bL \bW_2 +  \bS_{\bL}  \bXi_2 \label{eq:LLs22}
\end{eqnarray}
\end{proposition}
\textbf{Proof of Proposition 4.3}
By multiplying the first equation in (\ref{eq:reach}) with $\cO_1$ to the left we obtain:
$$
\mathcal{O}_1 \bA_1 \mathcal{R}_1  +  \mathcal{O}_1 \bK_2 \mathcal{R}_2 \bS_{\bR} + \mathcal{O}_1 \bB_1 \bR = \mathcal{O}_1 \bE_1 \mathcal{R}_1 \bLambda_1 \ \Rightarrow \ - \sIL_1 + \bXi_2 \bS_{\bR} +\bV_1 \bR = - \IL_1 \Lambda_1 
$$
and hence relation (\ref{eq:LLs11}) is proven. Similarly we prove (\ref{eq:LLs21}). By multiplying the first equation in (\ref{eq:obs}) with $\cR_1$ to the right we obtain:
$$
\mathcal{O}_1 \bA_1 \mathcal{R}_1  +  \bS_{\bL} \mathcal{O}_2 \bK_1 \mathcal{R}_1  + \bL \bC_1 \mathcal{R}_1  =\bM_1 \mathcal{O}_1 \bE_1 \mathcal{R}_1  \ \Rightarrow \ - \sIL_1 + \bS_{\bL} \bXi_1  +\bL \bW_1  = - \bM_1 \IL_1  
$$
and hence relation (\ref{eq:LLs12}) is proven. Similarly we prove (\ref{eq:LLs22}).
\noindent

\vspace{-2mm}

\begin{proposition}
The Loewner matrices $\IL_{1}$ and $\IL_{2}$ satisfy the following Sylvester equations:
\begin{eqnarray}
\bM_1 \IL_1 - \IL_1 \bLambda_1  &=& (\bV_1 \bR -\bL \bW_1)+( \bXi_2 \bS_\bR- \bS_\bL \bXi_1), \label{eqL1} \\[1mm]
\bM_2 \IL_2 - \IL_2 \bLambda_2  &=& (\bV_2 \bR - \bL \bW_2) + (\bXi_1 \bS_\bR - \bS_\bL \bXi_2) .
\label{eqL2}
\end{eqnarray} 
\end{proposition}
\textbf{Proof of Proposition 4.4}
By subtracting equation (\ref{eq:LLs11}) from (\ref{eq:LLs12}) we directly obtain (\ref{eqL1}) and also, by subtracting equation (\ref{eq:LLs21}) from (\ref{eq:LLs22})
we directly obtain (\ref{eqL2}).

\begin{proposition}
The shifted Loewner matrices $\sIL_{1}$ and $\sIL_{2}$ satisfy the following Sylvester equations:
\begin{eqnarray}
\bM_1 \sIL_1 - \sIL_1 \bLambda_1  &=& (\bM_1 \bV_1 \bR -\bL \bW_1 \bLambda_1)+(\bM_1 \bXi_2 \bS_\bR- \bS_\bL \bXi_1 \bLambda_1), \label{eqLs1} \\[1mm]
\bM_2 \sIL_2 - \sIL_2 \bLambda_2  &=& (\bM_2 \bV_2 \bR - \bL \bW_2 \bLambda_2) + (\bM_2 \bXi_1 \bS_\bR - \bS_\bL \bXi_2 \bLambda_2) .
\label{eqLs2}
\end{eqnarray} 
\end{proposition}
\textbf{Proof of Proposition 4.3}
By subtracting equation (\ref{eq:LLs11}) after being multiplied with $\bM_1$ to the left from equation (\ref{eq:LLs12}) after being multiplied with $\bLambda_1$ to the right, we directly obtain (\ref{eqLs1}). Similar procedure is applied to prove (\ref{eqLs2}).

\begin{remark}
{\rm The right hand side of the equations (\ref{eqL1}) - (\ref{eqLs2}) contains constant $0/1$ matrices  (i.e. $\bR,\bL,\bS_{\bR},\bS_{\bL}$) as well as matrices (i.e. $\bV_j,\bW_j,\bXi_j, \ j \in \{1,2\}$) which are directly constructed  by putting together the given samples values as pointed out in Example 4.2.}
\end{remark}

\subsection{Construction of reduced order models}

As we already noted, the interpolation data for the LSS case is significantly different than the one used for the linear case without switching, as higher order transfer function values are matched as shown in the previous sections. However, the rest of the procedure remains more or less unchanged.

\begin{lemma}
Assume that $k=\ell$ and that none of the
interpolation points $\lambda_i$, $\mu_j$ are eigenvalues of the pencils $\,(\sIL_1,\,\IL_1)$ and $\,(\sIL_2,\,\IL_2)$. Moreover,  consider the Loewner matrices $\IL_1$ and $\IL_2$ to be invertible. Then, a realization of a reduced order LSS  $\hat{\Si}$ that matches the data of the original LSS $\Si$ (as introduced in Definition 4.6) is given by the following matrices,
$$
\begin{cases} \hat{\bE}_1=-\IL_1, ~~\hat{\bA}_1=-\sIL_1, ~~\hat{\bB}_1=\bV_1, ~~\hat{\bC}_1=\bW_1, \\
\hat{\bE}_2=-\IL_2, ~~\hat{\bA}_2=-\sIL_2, ~~\hat{\bB}_2=\bV_2, ~~\hat{\bC}_2=\bW_2
\end{cases} \ \ \text{\rm and} \ \ \ \hat{\bK}_1=\bXi_1, \ \ \hat{\bK}_2=\bXi_2.
$$
If $k=n$, then the proposed realization is equivalent to the original one (as in Definition 4.1). 
\end{lemma}
\textbf{Proof of Lemma 4.4} This result directly follows from Lemma 4.1 by taking into consideration the notations introduced in (\ref{loewner}-\ref{loewner2}).\\

\vspace{-2mm}

In the case of redundant data, at least one of the pencils $(\sIL_j,\,\IL_j)$ is singular (for $j \in \{1,2\}$), and  hence construct pairs of projectors $(\bX_j,\bY_j)$ (corresponding to mode $j$) similar to (\ref{projection}). The MOR procedure for approximate data matching is presented as follows.

\vspace{2mm}

\noindent
{\bf{Procedure 1}} \ Consider the rank revealing singular value factorization of the  Loewner matrices, 
\begin{equation}
\IL_j  = \left[ \begin{array}{cc}
\bY_j & \tilde{\bY}_j \end{array} \right]  \left[ \begin{array}{cc}
\bS_j & \bO \\ \bO & \tilde{\bS}_j \end{array} \right] \left[ \begin{array}{cc}
\bX_j & \tilde{\bX}_j \end{array} \right]^T =  \bY_j \bS_j \bX_j^T +\tilde{\bY}_j \tilde{\bS}_j \tilde{\bX}_j^T ,
\end{equation}
where $\bY_j, \bX_j \in \mathbb{R}^{k \times r_j}$ and $\bS_j \in \mathbb{R}^{r_j \times r_j}$.
The projected system matrices corresponding to subsystem $\hat{\Si}_j$ are computed as,
\vspace{-1mm}
$$\hat{\bE}_j = - \bY_j^{T} \IL_j \bX_j, ~\hat{\bA}_j = -\bY_j^{T} \sIL_j \bX_j, ~\hat{\bB}_j = \bY_j^{T} \bV_j,~ 
\hat{\bC}_j = \bW_j \bX_j, \ \ \text{for} \ j \in \{1,2\}
\vspace{-1mm}
$$ 
Moreover, the projected coupling matrices  are computed in the following way
$$
\hat{\bK}_1 = \bY_2^{T} \bXi_1 \bX_1, \ \ \hat{\bK}_2 = \bY_1^{T} \bXi_2 \bX_2.
$$
By choosing $r_j$ as the numerical rank of the Loewner matrix $\IL_j$ (i.e. the largest neglected singular value corresponding to index $r_j+1$ is less than machine precision $\epsilon$), ensure that the $\hat{\bE}_j$ matrices are not singular. Hence, construct a reduced order LSS denoted with $\hat{\Si}$, that approximately matches the data of the original LSS $\Si$.
If the truncated singular values are all 0 (the ones on the main diagonal of the matrices $\tilde{\bS}_j$), then the matching is exact.

We provide a qualitative rather than quantitative result for the projected Loewner case. The quality of approximation is directly linked to the singular values of the Loewner pencils which represent an indicator of the desired accuracy. For linear systems with no switching, an error bound is provided in \cite{birkjour} as a quantitative measure.

The dimensions of the subsystems $\hat{\Si}_1$ and $\hat{\Si}_2$, corresponding to the reduced order LSS, need not be the same (i.e. $r_1 \neq r_2$). In this case the coupling matrices are not square anymore.

The projectors are computed via singular value factorization of the Loewner matrices.
The use of the \textit{Drazin} or \textit{Moore-Penrose} pseudo inverses also holds (as shown in \cite{drazin}).

\section{The Loewner framework for linear switched systems - the general case} \label{sec:main}

In this section we are mainly concerned with generalizing some of the results presented in Section 4. Most of the findings can be smoothly extended to the cases with more complex switching patterns (more modes). The key for this is enforcing a cyclic structure of the interpolation framework, so that, everything can be written in matrix equation format.
 
\begin{definition}
Let $\Gamma$ and $\Theta$ be finite sets of tuples so that $\Gamma , \Theta \subseteq \displaystyle \bigcup_{k=1}^{\infty} Q^k \times \mathbb{C}^k$ so that $\Gamma$ has the prefix closure property, i.e.
 $$
 (q_1,q_2, \ldots,q_i,\lambda_1,  \ldots,\lambda_i) \in \Gamma \Rightarrow (q_2, \ldots,q_i,\lambda_2, \ldots,\lambda_i) \in \Gamma \ \forall i \geqslant 2
 $$
and $\Theta$ has the suffix closure property, i.e.
 $$
 (q_1,q_2, \ldots,q_j,\mu_1, \ldots,\mu_j) \in \Theta \Rightarrow (q_1, \ldots,q_{j-1},\mu_1, \ldots,\mu_{j-1}) \in \Theta \ \forall j \geqslant 2
 $$
\end{definition}
\vspace{-2mm}
\noindent
Now consider the specific subset $\Gamma_q$ (for any $q \in Q$) of the set $\Gamma$ in the following way:
\vspace{-1mm}
$$
\Gamma_q = \{ (q_1,q_2, \ldots,q_i,\lambda_1, \ldots,\lambda_i) \in \Gamma \ \vert \ q_1 = q, \ i \leqslant \delta_{\Gamma} \}, \ \ \ \delta_{\Gamma} = \underset{\ w \in \Gamma}{\max(\vert w \vert)}/2
\vspace{-3mm}
$$
 Denote the cardinality of $\Gamma_q$ with $k_q = \text{card}(\Gamma_q)$ and explicitly enumerate the elements of this set as follows: $\Gamma_q = \{ w_q^{(1)},w_q^{(2)}, \ldots  ,w_q^{(k_q)}\}$. Consider the following function (mapping) $\br:\Gamma_q \rightarrow \mathbb{C}^{n_q \times 1}$ that maps a word form $\Gamma_q$ into a column vector of size $n_q$:
\vspace{-1mm}
$$
\br((q,q_2,\ldots,q_i,\lambda_1,\ldots,\lambda_i)) = \bPhi_{q}(\lambda_1) \bK_{q_2,q} \bPhi_{q_2}(\lambda_2) \cdots \bK_{q_{i},q_{i-1}}\bPhi_{q_i}(\lambda_i) \bB_{q_i}
\vspace{-1mm}
$$
Now we are ready to construct the reachability matrix $\cR_q$ corresponding to the mode $q$ of the system $\Si$ as follows:
\begin{equation}\label{eq:R_def_gen}
\cR_q = \left[ \begin{array}{cccc}
\br(w_q^{(1)}) & \br(w_q^{(2)})  & \cdots & \br(w_q^{(k_q)}) 
\end{array}  \right] \in \mathbb{C}^{n_q \times k_q}
\end{equation}

\vspace{-2mm}
\noindent
Similarly, define the specific subset $\Theta_q$ (for any $q \in Q$) of the set $\Theta$ in the following way:
\vspace{-1mm}
$$
\Theta_q = \{ (q_1,q_2,\ldots,q_j,\mu_1,\ldots,\mu_j) \in \Gamma \ \vert \ q_j = q, \ j \leqslant \delta_{\Theta} \},  \ \ \delta_{\Theta} = \underset{\ w \in \Theta}{\max(\vert w \vert)}/2
\vspace{-2mm}
$$
Consider the cardinality of $\Theta_q$ to be the same as the one of $\Gamma_q$, i.e. $k_q = \text{card}(\Theta_q)$. Although this additional constraint is not necessarily needed, we would like to enforce the construction of reduced systems with square matrices $\bA_k$ and $\bE_k$. Next we explicitly enumerate the elements of this set as follows: $\Theta_q = \{ v_q^{(1)},v_q^{(2)},\ \ldots ,v_q^{(k_q)}\}$. Consider the following mapping $\boo:\Theta_q \rightarrow \mathbb{C}^{1 \times n_q}$ that maps a word form $\Theta_q$ into a row vector of size $n_q$:
\vspace{-1mm}
$$
\boo((q_1,q_2,\ldots,q_{j-1},q,\mu_1,\ldots,\mu_j)) =\bC_{q_1} \bPhi_{q_1}(\mu_1) \bK_{q_2,q_1} \bPhi_{q_2}(\mu_2) \cdots \bK_{q,q_{j-1}}\bPhi_{q}(\mu_j) 
\vspace{-1mm}
$$
Now we are ready to construct the observability matrix $\cO_q \in \mathbb{C}^{k_q \times n_q}$ corresponding to the mode $q$ of the system $\Si$ as follows
\begin{equation}\label{eq:O_def_gen}
\cO_q = \left[ \begin{array}{cccc}
{\boo(v_q^{(1)})}^{T} & {\boo(v_q^{(2)})}^{T}   & \cdots & {\boo(v_q^{(k_q)})}^{T} 
\end{array}  \right]^T \in \mathbb{C}^{ k_q \times n_q}
\vspace{-1mm}
\end{equation}
Consider the following example to show how the general procedure is extended from the linear case (no switching) to the case when switching occurs.

\begin{example}
Take $D=3$ (3 active modes) and hence $Q = \{1,2,3\}$. The following interpolation points are given: $\{s_1,s_2,\ldots,s_{18}\} \subset \mathbb{C}$. The first step is to partition this set into two disjoint subsets (each having 9 points):
$$
 \text{\underline{left interpolation points}} : \{\mu_1,\mu_2,\ldots,\mu_9\}  \ \ \ \ \ \ \ \text{\underline{right interpolation points}} : \{\lambda_1,\lambda_2,\ldots,\lambda_9\} 
$$
The set $\Gamma$ is composed of three subsets $\Gamma = \Gamma_1 \bigcup \Gamma_2 \bigcup \Gamma_3$ which are defined in a cyclic way by imposing the previously defined suffix closure property, as follows
 $$\begin{cases} \Gamma_1 = \{ (1,\lambda_1), (1,3,\lambda_4,\lambda_3), (1,3,2,\lambda_7,\lambda_6,\lambda_2) \} \\ \Gamma_2 = \{ (2,\lambda_2), (2,1,\lambda_5,\lambda_1), (2,1,3,\lambda_8,\lambda_4,\lambda_3) \} \\ \Gamma_3 = \{ (3,\lambda_3), (3,2,\lambda_6,\lambda_2), (3,2,1,\lambda_9,\lambda_5,\lambda_1) \} \end{cases}
 $$ 
  To the sets $\Gamma_j$, we attach the reachability matrices $\cR_j$ defined as follows:
\begin{eqnarray*}
&& \cR_1 = \left[ \begin{array}{ccc} \bPhi_1(\lambda_1) &  \bPhi_1(\lambda_4) \bK_{3,1} \bPhi_3(\lambda_3) \bB_3 &  \bPhi_1(\lambda_7) \bK_{3,1} \bPhi_3(\lambda_6) \bK_{2,3} \bPhi_2(\lambda_2) \bB_2
\end{array} \right], \\
&& \cR_2 = \left[ \begin{array}{ccc} \bPhi_2(\lambda_2) &  \bPhi_2(\lambda_5) \bK_{1,2} \bPhi_1(\lambda_1) \bB_1 &  \bPhi_2(\lambda_8) \bK_{1,2} \bPhi_1(\lambda_4) \bK_{3,1} \bPhi_3(\lambda_3) \bB_3
\end{array} \right], \\
&& \cR_3 = \left[ \begin{array}{ccc} \bPhi_3(\lambda_3) &  \bPhi_3(\lambda_6) \bK_{2,3} \bPhi_2(\lambda_2) \bB_2 &  \bPhi_3(\lambda_9) \bK_{2,3} \bPhi_2(\lambda_5) \bK_{1,2} \bPhi_1(\lambda_1) \bB_1
\end{array} \right]
\end{eqnarray*}
In the same manner, the set $\Theta$ is composed of three subsets $\Theta = \Theta_1 \bigcup \Theta_2 \bigcup \Theta_3$ which are again defined in a cyclic way by imposing the previously defined prefix closure property, as follows
 $$
 \begin{cases} \Theta_1 = \{ (1,\mu_1), (3,1,\mu_3,\mu_4), (1,2,1,\mu_1,\mu_5,\mu_7) \} \\ \Theta_2 = \{ (2,\mu_2), (1,2,\mu_1,\mu_5), (2,3,2,\mu_2,\mu_6,\mu_8) \} \\ \Theta_3 = \{ (3,\mu_3), (2,3,\mu_2,\mu_6), (3,1,3,\mu_3,\mu_4,\mu_9) \} \end{cases}
 $$
 To the sets $\Theta_i$, we attach the observability matrices $\cO_i$ defined as follows:
$$
\cO_1 = \left[ \begin{array}{c} \bC_1 \bPhi_1(\mu_1) \\ \bC_3 \bPhi_3(\mu_3) \bK_{1,3} \bPhi_1(\mu_4) \\ \bC_1 \bPhi_1(\mu_1) \bK_{2,1} \bPhi_2(\mu_5) \bK_{1,2} \bPhi_1(\mu_7)
\end{array} \right], \ \ \cO_2 = \left[ \begin{array}{c} \bC_2 \bPhi_2(\mu_2) \\ \bC_1 \bPhi_1(\mu_1) \bK_{2,1} \bPhi_2(\mu_5) \\ \bC_2 \bPhi_2(\mu_2) \bK_{3,2} \bPhi_3(\mu_6) \bK_{2,3} \bPhi_2(\mu_8)
\end{array} \right]
$$

\newpage

$$
\cO_3 = \left[ \begin{array}{c} \bC_3 \bPhi_3(\mu_3) \\ \bC_2 \bPhi_2(\mu_2) \bK_{3,2} \bPhi_3(\mu_6) \\ \bC_3 \bPhi_3(\mu_3) \bK_{1,3} \bPhi_1(\mu_4) \bK_{3,1} \bPhi_3(\mu_9)
\end{array} \right]
$$
\end{example}


\subsection{Sylvester equations for $\cR_q$ and $\cO_q$}

In this section we would like to generalize the results presented in Lemma 4.2 and Lemma 4.3, and hence extend the framework to a general number of operational modes denoted with D.

\begin{definition}
Introduce the special concatenation of tuples composed of mixed elements (symbols) that are from two different sets (for example $Q$ and $\mathbb{C}$) as the mapping  with the following property: 
$$
\big{(} \alpha_1 \circledcirc \beta_1  \big{)}\odot \big{(} \alpha_2 \circledcirc \beta_2 \big{)} = \Big{(} \big{(} \alpha_1  \circledcirc \alpha_2 \big{)} \circledcirc  (\beta_1 \circledcirc \beta_2 \big{)} \Big{)}, \ \ 
$$
where $\alpha_k \in Q^{i_k}$ and $\mathbf{\beta}_k \in \mathbb{C}^{j_k}$ for $k = 1,2$.

\end{definition}

\begin{definition}
For $g,i = 1,\ldots,D$, let $\bS_i^{(g)} = \left[ \begin{array}{ccc}
\bS_i^{(g)}(1) & \ldots & \bS_i^{(g)}(k_g)
\end{array} \right] \in \mathbb{R}^{k_i \times k_g}$ be constant matrices that contain only $0 / 1$ entries constructed so that $ \bS_i^{(g)}(1) = \bfz_{k_i} $ and for $u = 2,\ldots,k_g$, we write:  
\begin{equation}\label{def_S_mat}
\bS_i^{(g)}(u) = \begin{cases} \bee_{u-1,k_i}, \ \text{\rm if } \ \exists \ \tilde{\lambda} \in \mathbb{C}, \ \text{\rm s.t.}  \ \bw_g^{(u)} = (g,\tilde{\lambda}) \odot \bw_{i}^{(u-1)}, \\
 \bfz_{k_i}, \ \text{\rm else} \end{cases}
\end{equation}
\end{definition}
Also, introduce the matrices $\bR^{(i)}$ and $\bLambda_i$ that are defined similarly as in (\ref{RMat}) and (\ref{LLam}), i.e.,
\begin{equation} \label{RLMat_gen}
\bR^{(i)}=\left[\bfe^T_{1,m_1}~\cdots~\bfe^T_{1,m_{k^\dagger}}\right] \in\IR^{1 \times k_i},\ \ \ \bLambda_i = \mbox{blkdiag}\,[\bLambda_i^{(1)},~~  \bLambda_i^{(2)},~\ldots,~ \bLambda_i^{(k^\dagger)}] \in\IR^{k_i \times k_i},
\end{equation}
where the diagonal matrices $\bLambda_i^{(a)}$, $a = 1,\ldots,k^{\dagger}$ contain the right interpolation points associated to mode $i$. For general cyclic structure incorporated of the set $\Gamma$, it follows that the reachability matrices $\cR_i \in \mathbb{R}^{n_i \times k_i}, \ 1 \leqslant i \leqslant D$ satisfy the following system of generalized Sylvester equations:
\begin{equation}
\begin{cases} 
 \bA_1 \cR_1 + \displaystyle \sum_{i=1}^{D} \bK_{i,1}  \cR_i \bS_i^{(1)}  +  \bB_1 \bR^{(1)} = \bE_1 \cR_1 \bLambda_1 \\ 
\bA_2 \cR_2  + \displaystyle \sum_{i=1}^{D} \bK_{i,2}  \cR_i \bS_i^{(2)}  + \bB_2 \bR^{(2)} = \bE_2 \cR_2 \bLambda_2 \\
\ \ \ \ \ \ \ \ \ \ \vdots \\
\bA_D   \cR_D + \displaystyle \sum_{i=1}^{D}  \bK_{i,D} \cR_i \bS_i^{(D)}  +  \bB_D \bR^{(D)} = \bE_D \cR_D \bLambda_D
\end{cases}
\end{equation}

\noindent
Note that $\bS_i^{(i)} = \bfz_{k_i,k_i}$, and if $k_1 = k_2 = \cdots = k_D = k$, the above defined matrices $\bS_i^{(g)}$ satisfy the following equality   $\forall g \in Q$ :
\begin{equation} \label{S_cond}
\displaystyle \sum_{i=1}^D \bS_i^{(g)} = ~\mbox{blkdiag}\,\left[\bJ_{m_1},~\ldots,\bJ_{m_{k^\dagger}}\right]
\vspace{-1mm}
\end{equation}
where $\bJ_{l}$ is the Jordan block of size $l$ defined in (\ref{SMat}).

To directly find $\cR_{g}, \ g =1,2,3$ for the case presented in Example 5.1, we have to solve the following system of coupled generalized Sylvester equations
$$
\begin{cases} \bA_1 \cR_1+\bK_{3,1} \cR_3 \bS_3^{(1)} + \bB_1 \bR = \bE_1 \cR_1 \bLambda_1 \\ \bA_2 \cR_2+\bK_{1,2} \cR_1 \bS_1^{(2)} + \bB_2 \bR = \bE_2 \cR_2 \bLambda_2 \\ \bA_3 \cR_3+\bK_{2,3} \cR_2 \bS_2^{(3)} + \bB_3 \bR = \bE_3 \cR_3 \bLambda_3
\end{cases}
$$
where:
$$
\bLambda_1 = \left[ \begin{array}{ccc} \lambda_1 & 0 & 0 \\ 0 & \lambda_4 & 0 \\ 0 & 0 & \lambda_7  \end{array} \right], \  \bLambda_2 = \left[ \begin{array}{ccc} \lambda_2 & 0 & 0 \\ 0 & \lambda_5 & 0 \\ 0 & 0 & \lambda_8  \end{array} \right], \  \bLambda_3 = \left[ \begin{array}{ccc} \lambda_3 & 0 & 0 \\ 0 & \lambda_6 & 0 \\ 0 & 0 & \lambda_9  \end{array} \right], 
$$
$$
\bR  =  \left[ \begin{array}{ccc} 1 & 0 & 0 \end{array} \right], \ \bS_3^{(1)} = \bS_1^{(2)} = \bS_2^{(3)} = \left[ \begin{array}{ccc} 0 & 1 & 0 \\ 0 & 0 & 1 \\ 0 & 0 & 0  \end{array} \right],
$$
This corresponds to the case $k_1 = k_2 = k_3=3 $, $k^{\dagger} = 1$ and $m_1 = 3$.

\begin{definition}
For $h,j = 1,\ldots,D$, let $\bT_j^{(h)} = \left[ \begin{array}{ccc}
\big{(} \bT_j^{(h)} \big{)}^T(1) & \ldots & \big{(} \bT_j^{(h)} \big{)}^T(k_h)
\end{array} \right]^T \in \mathbb{R}^{\ell_h \times \ell_j}$ be constant matrices that contain only $0 / 1$ entries constructed so that $ \big{(} \bT_j^{(h)} \big{)}^T(1) = \bfz_{\ell_j} $ and for $v = 2,\ldots,k_g$, we write:  
\begin{equation}\label{def_T_mat}
\big{(} \bT_j^{(h)} \big{)}^T(v) = \begin{cases} \bee_{v-1,k_j}, \ \text{\rm if } \ \exists \ \tilde{\mu} \in \mathbb{C}, \ \text{\rm s.t.}  \ \bw_h^{(v)} =   \bw_{j}^{(v-1)} \odot (h,\tilde{\mu}), \\
 \bfz_{\ell_j}, \ \text{\rm else} \end{cases}
\end{equation}
\end{definition}
\normalsize
Also, introduce the following matrices
\begin{equation} \label{MLam_gen}
\big{(} \bL^{(j)} \big{)}^T=\left[\bfe^T_{1,p_1}~\cdots~\bfe^T_{1,p_{\ell^\dagger}}\right] \in\IR^{1 \times \ell_j}, \ \ \bM_j = \mbox{blkdiag}\,[\bM_j^{(1)},~~\bM_j^{(2)},~\ldots,~\bM_j^{(\ell^\dagger)}] \in  \IR^{\ell_j \times \ell_j},
\end{equation} 
where the diagonal matrices $\bM_j^{(v)}$ for $v = 1,\ldots,\ell_j$ contain the left interpolation points associated to mode $j$. For general cyclic structure incorporated by definition in the set $\Theta$, one can conclude that the observability matrices $\cO_j \in \mathbb{R}^{\ell_j \times n_j}, \ 1 \leqslant j \leqslant D$ satisfy the following system of generalized Sylvester equations:
\begin{equation}
\begin{cases} 
\cO_1 \bA_1 + \displaystyle \sum_{j=1}^{D} \bT_j^{(1)} \cO_j \bK_{1,j} + \bL^{(1)} \bC_1 = \bM_1 \cO_1 \bE_1 \\ 
\cO_2 \bA_2 + \displaystyle \sum_{j=1}^{D} \bT_j^{(2)} \cO_j \bK_{2,j} + \bL^{(2)} \bC_2 = \bM_2 \cO_2 \bE_2 \\
\ \ \ \ \ \ \ \ \ \ \vdots \\
\cO_D \bA_D + \displaystyle \sum_{j=1}^{D} \bT_j^{(D)} \cO_j \bK_{D,j} + \bL^{(D)} \bC_D = \bM_D \cO_D \bE_D
\end{cases}
\end{equation}
Note that $\bT_j^{(j)} = \bfz_{\ell_j,\ell_j}$, and if $\ell_1 = \ell_2 = \cdots = \ell_D = \ell$, the square matrices $\bT_j^{(h)} \in \mathbb{R}^{\ell \times \ell}$ satisfy the following equality,  $\forall h \in Q$:
\begin{equation} \label{T_cond}
\displaystyle \sum_{j=1}^D \bT_j^{(h)} = ~\mbox{blkdiag}\,\left[\bJ_{p_1},~\ldots,\bJ_{p_{l^\dagger}}\right]^T
\end{equation} 

Again to find the matrices $\cO_h, \ h =1,2,3$ in Example 5.1, it is required to solve the following system of coupled generalized Sylvester equations
$$
\begin{cases}
\cO_1 \bA_1+ \bT_3^{(1)} \cO_3 \bK_{1,3} + \bT_2^{(1)} \cO_2 \bK_{1,2}  + \bL \bC_1 = \bM_1 \cO_1 \bE_1 \\ \cO_2 \bA_2+ \bT_1^{(2)} \cO_1 \bK_{2,1} + \bT_3^{(2)} \cO_3 \bK_{2,3}  + \bL \bC_2 = \bM_2 \cO_2 \bE_2 \\ \cO_3 \bA_3+ \bT_2^{(3)} \cO_2 \bK_{3,2} + \bT_1^{(3)} \cO_1 \bK_{3,1}  + \bL \bC_3 = \bM_3 \cO_3 \bE_3
\end{cases}
$$
where:
$$
\bM_1 = \left[ \begin{array}{ccc} \mu_1 & 0 & 0 \\ 0 & \mu_4 & 0 \\ 0 & 0 & \mu_7  \end{array} \right], \  \bM_2 = \left[ \begin{array}{ccc} \mu_2 & 0 & 0 \\ 0 & \mu_5 & 0 \\ 0 & 0 & \mu_8  \end{array} \right], \  \bM_3 = \left[ \begin{array}{ccc} \mu_3 & 0 & 0 \\ 0 & \mu_6 & 0 \\ 0 & 0 & \mu_9  \end{array} \right], 
$$
$$
 \bT_3^{(1)} = \bT_1^{(2)} = \bT_2^{(3)} = \left[ \begin{array}{ccc} 0 & 0 & 0 \\ 1 & 0 & 0 \\ 0 & 0 & 0  \end{array} \right], \    \bT_2^{(1)} = \bT_3^{(2)} = \bT_1^{(3)} = \left[ \begin{array}{ccc} 0 & 0 & 0 \\ 0 & 0 & 0 \\ 0 & 1 & 0  \end{array} \right], \ \bL =  \bfe_1
$$
This corresponds to the case $\ell_1 = \ell_2 = \ell_3 = 3$, $l^{\dagger} = 1$ and $p_1 = 3$. Note that the relation in (\ref{T_cond}) \\[1mm] hold, i.e., $
 \bT_2^{(1)}+ \bT_3^{(1)} =  \bT_1^{(2)} + \bT_3^{(2)} =  \bT_1^{(3)} +  \bT_2^{(3)} = \bJ_3^T$.

\subsection{The Loewner matrices}

For the case of linear switched systems with D active modes, the generalization of the Loewner framework includes one important feature. Instead of only one pair of Loewner matrices (as in the linear case without switching which is covered in Section 3), we define a pair of Loewner matrices for each individual active mode; hence in total D pairs of Loewner matrices.

\begin{definition} Given a linear switched system $\Si$, let $\{\cR_i \vert i \in Q\}$ and $\{\cO_j \vert j \in Q\}$ be the controllability and observability matrices associated with the multi-tuples $\Gamma_i$ and $\Theta_j$. The {\bf Loewner matrices} $\{ \IL_i \vert \ i \in Q \}$ are defined as
\begin{equation} \label{eq:Loew_gen}
\IL_1= - \cO_1\,\bE_1\,\cR_1,~\IL_2= - \cO_2\,\bE_2\,\cR_2, \ \ldots , \ \IL_D= - \cO_D\,\bE_D\,\cR_D
\end{equation}
Additionally, the {\bf shifted Loewner matrices} $\{ \sIL_i \vert \ i \in Q \}$ are defined as
\begin{equation} \label{eq:SLoew_gen}
\sIL_1= - \cO_1\,\bA_1\,\cR_1,~~\sIL_2= - \cO_2\,\bA_2\,\cR_2, \ \ldots , \ \sIL_D= - \cO_D\,\bA_D\,\cR_D
\vspace{-2mm}
\end{equation}
Also introduce the matrices $ \ \forall i,j \in Q$
\vspace{-2mm}
\begin{equation*}
\bW_i=\bC_i\,\cR_i,~~ \bV_j= \cO_j\,\bB_j, ~~\mbox{and}~~  \bXi_{i,j}=\cO_j\,\bK_{i,j}\, \cR_i 
\end{equation*}
\end{definition}

\begin{remark}
{\rm The number of Loewner matrices, shifted Loewner matrices, $W_i$ row vectors and $V_j$ column vectors is the same as the number of active modes (i.e D). On the other hand, the number of matrices $\bXi_{i,j}$ increases quadratically with D (i.e in total $D^2$ matrices).}
\end{remark}

\begin{remark}
{\rm Note that the matrices $\IL_i$ and $\sIL_i$ as defined in (\ref{eq:Loew_gen}) and (\ref{eq:SLoew_gen}) (for $i \in \{1,2,\ldots,D\}$) are indeed Loewner matrices, that is, they can be expressed as {\it divided differences} of generalized transfer function values of the underlying LSS.}
\end{remark}

\begin{proposition}
The Loewner matrices $\IL_{h}$ satisfy the following Sylvester equations:
\begin{equation}\label{eqLg}
\bM_{h} \IL_{h} - \IL_{h} \bLambda_{h}  = (\bV_{h} \bR-\bL \bW_{h}) + \displaystyle \sum_{j=1}^{D} \Big{(}  \bXi_{j,h} \bS_j^{(h)} -  \bT_j^{(h)} \bXi_{h,j} \Big{)}, \ h \in Q
\end{equation} 
\end{proposition}

\begin{proposition}
The shifted Loewner matrices $\sIL_{h}$ satisfy the following Sylvester equations:
\begin{equation}\label{eqLsg}
\bM_{h} \sIL_{h} - \sIL_{h} \bLambda_{h} = (\bM_{h} \bV_{h} \bR - \bL \bW_{h} \bLambda_{h}) + \displaystyle \sum_{j=1}^{D} \Big{(} \bM_{h} \bXi_{j,h} \bS_j^{(h)} - \bT_j^{(h)} \bXi_{h,j} \bLambda_{h}  \Big{)}, \ h \in Q
\end{equation} 
\end{proposition}

\begin{remark}
{\rm The proof of the results stated in (\ref{eqLg})-(\ref{eqLsg}) is performed in a similar manner as for the results obtained for the special case $D=2$ in Section 4 (i.e. for (\ref{eqL1})-(\ref{eqLs2})). }
\end{remark} 
 
 \subsection{Construction of reduced order models}

The general procedure for the case with $D$ switching modes is more or less similar to the one covered in Section 4.3 (where $D=2$).  \\[-4mm]

\begin{lemma}
 Let $\IL_j$ be invertible matrices for $j=1,\ldots,D$, such that none of the interpolation points $\lambda_i$, $\mu_k$ are eigenvalues of any of the Loewner pencils $(\sIL_j,\,\IL_j)$. Then, it follows that the matrices
$$
\{ \hat{\bE}_j=-\IL_j, ~~\hat{\bA}_1=-\sIL_j, ~~\hat{\bB}_j=\bV_j, ~~\hat{\bC}_j=\bW_j, ~~ \hat{\bK}_{i,j}=\bXi_{i,j} \}, \ \ i,j \in \{1,\ldots,D\}
 \vspace{-1mm}
$$
form a realization of a reduced order LSS $\hat{\Si}$ that matches the data of the original LSS $\Si$. 
 If $k_j=n_j$ for $j = 1,\ldots,D$, the proposed realization is equivalent to the original one. 
\end{lemma}
The concept of a LSS matching the data of another LSS in the case $D>2$ is formulated in a similar manner as to the case $D=2$ (i.e., which is covered in Definition 4.6). Also, the definition of equivalent LSS for the case $D>2$ is formulated similarly as to Definition 4.1.

In the case of redundant data, at least one of the pencils $(\sIL_j,\,\IL_j)$ is singular (for $j \in \{1,\ldots,D\}$). The main procedure is presented as follows.

\vspace{2mm}

\noindent
{\bf{Procedure 2}} \ Consider the rank revealing singular value factorization of the Loewner matrices:
\begin{equation}
\IL_j  = \left[ \begin{array}{cc}
\bY_j & \tilde{\bY}_j \end{array} \right]  \left[ \begin{array}{cc}
\bS_j & \bO \\ \bO & \tilde{\bS}_j \end{array} \right] \left[ \begin{array}{cc}
\bX_j & \tilde{\bX}_j \end{array} \right]^T =  \bY_j \Si_j \bX_j^T + \tilde{\bY}_j \tilde{\bS}_j \tilde{\bX}_j^T
\end{equation}
where $\bX_j,\bY_j \in \mathbb{R} ^{k_j \times r_j}$ , $\bS_j \in \mathbb{R} ^{r_j \times r_j}$ and $j = \{1,\ldots,D\}$. Here, choose $r_j$ as the numerical rank of the Loewner matrix $\IL_j$ (i.e. the largest neglected singular value corresponding to index $r_j+1$ is less than machine precision $\epsilon$).
The projected system matrices computed as 
\vspace{-1mm}
$$
\hat{\bE}_j = - \bY_j^{T} \IL_j \bX_j, ~\hat{\bA}_j = -\bY_j^{T} \sIL_j \bX_j, ~\hat{\bB}_j = \bY_j^{T} \bV_j,~ 
\hat{\bC}_j = \bW_j \bX_j, \ \ \text{for} \ j \in \{1,\ldots,D\}
\vspace{-1mm}
$$ 
and the projected coupling matrices computed in the following way
$$
\hat{\bK}_{i,j} = \bY_j^{T} \bXi_{i,j} \bX_i, \ \ \forall i,j \in \{1,\ldots,D\},
$$ 
form a realization of a reduced order LSS denoted with $\hat{\Si}$ that approximately matches the data of the original LSS $\Si$. Each reduced subsystem $\hat{\Si}_j$ has dimension $r_j, \ j \in \{1,\ldots,D\}$.

\begin{remark}
{\rm 
If the truncated singular values are all 0 (the ones on the main diagonal of the matrices $\tilde{\bS}_j$), then the interpolation is exact.}
\end{remark}


\newpage

\section{Numerical experiments} \label{sec:num}
In this section we illustrate the new method by means of three numerical examples. We use a certain generalization of the balanced truncation (BT) method for LSS (as presented in (\cite{mtc12}) to compare the performance of our new introduced method. The main ingredient of the BT method is to compute the
 the controllability and observability gramians $\bP_i$ and $\bQ_i$ (where $i \in \{1,2,\ldots,D\}$) as the solutions of the following Lyapunov equations:
 \begin{eqnarray}
&& \bA_i \bP_i \bE_i^T + \bE_i \bP_i \bA_i^T + \bB_i \bB_i^T = \b0 \\
&& \bA_i^T \bQ_i \bE_i + \bE_i^T \bP_i \bA_i + \bC_i^T \bC_i = \b0 
 \end{eqnarray}
 
\subsection{Balanced Truncation} 
In \cite{mtc12} it has been shown that, if certain conditions
are satisfied, the technique of simultaneous balanced truncation
can be applied to switched linear systems. Hence, in some special cases, the existence of a global transformation matrix $\bV_{bal}$ is guaranteed; it follows that:
\begin{equation}
\bV_{bal} \bP \bV_{bal}^T = \bV_{bal}^{-T} \bQ V_{bal}^{-1} = \bU_i
\end{equation}
where $\bU_i$ are diagonal matrices. Although conceptually attractive as a MOR method, in general the conditions are rather restrictive in practice. This motivates the search for a more
general MOR approach for the case where simultaneous
balancing cannot be achieved.

The problem of finding a balancing transformation for a single linear system can be formulated as finding a nonsingular matrix such that the following cost function is minimized (see \cite{book}):
\begin{equation}
f(\bV) = \text{trace} [ \bV \bP \bV^T + \bV^{-T} \bQ \bV^{-1}] 
\end{equation}
For the class of LSS with distinct operational modes, we hence have to minimize not one but a number of D cost functions:
\begin{equation}\label{eq:f_dist_modes}
f_i(\bV) = \text{trace} [ \bV \bP_i \bV^T + \bV^{-T} \bQ_i \bV^{-1}] , \ \ i \in \{1,2,\ldots,D\}
\end{equation}
If the conditions of Corollary IV.3 from \cite{mtc12} hold,
simultaneous balancing is possible, and there exists a transformation $\bV$ which simultaneously minimizes $f_i$ for all $i = 1,2, \ldots,D$. Instead of having D functions as in (\ref{eq:f_dist_modes}), one can introduce a single overall cost function (i.e the average of the cost functions of the individual modes). Define the function $f_{av}$ as in \cite{mtc12}:
\begin{equation}\label{eq:cost_fct}
f_{av}(\bV) = \frac{1}{D} \sum_{i=1}^D	\text{trace} [ \bV \bP_i \bV^T + \bV^{-T} \bQ_i \bV^{-1}]  = \text{trace} [ \bV \bP_{av} \bV^T + \bV^{-T} \bQ_{av} \bV^{-1}] 
\end{equation}
where
\begin{equation}
\bP_{av} = \frac{1}{D} \sum_{i=1}^D \bP_i, \ \ \bQ_{av} = \frac{1}{D} \sum_{i=1}^D \bQ_i,
\end{equation}
In the case of LSS, the BT method computes a basis where the sum of the sum of the eigenvalues of $\bP_i$ and $\bQ_i$ over all modes is minimal. Hence, minimizing the proposed overall cost function provides a natural extension of classical BT to the case of LSS.

It follows that the transformation $\tilde{\bV}$ that minimizes the cost function in (\ref{eq:cost_fct}) is precisely the one which balances the pair $(\bP_{av},\bQ_{av})$ of average gramians.

By applying $\tilde{\bV}$ to the individual modes and truncating, a reduced order model is obtained. After applying the transformation $\tilde{\bV}$, the new state space representations of the individual modes need not be balanced. Nevertheless, as stated in  \cite{mtc12}, it is expected to be relatively close to being balanced.	


\subsection{First example}
As first example we consider the simple model of an evaporator vessel from \cite{mzb99}. There is a
constant inflow of liquid fin into a tank and an outflow that depends on the pressure in the
tank and the Bernoulli resistance $R_b$. To keep the level of fluid in the evaporator vessel at or
below a pre-specified maximum, an overflow mechanism is activated when the level of fluid L in
the evaporator exceeds the threshold value $L_{th}$. This causes a flow through a narrow pipe with resistance $R_p$ and inertia I that builds up flow momentum p. The system is modeled in two distinct operation modes: mode 1, where there is no overflow (the fluid level is below the overflow level), and mode 2, where the overflow mechanism is active. The ordinary differential equations describing the system in the two operation modes are given by
\begin{eqnarray*}
&& \left[ \begin{array}{cc}
I & 0 \\ 0 & C
\end{array} \right] \left[ \begin{array}{cc}
\dot{p} \\ \dot{L}
\end{array} \right] = \left[ \begin{array}{cc}
-R_p & 0 \\ 0 & -1/R_b
\end{array} \right] \left[ \begin{array}{cc}
\dot{p} \\ \dot{L}
\end{array} \right] + \left[ \begin{array}{cc}
0 \\ f_{in}
\end{array} \right] \ \ (\text{mode} \ 1) \\
&& \left[ \begin{array}{cc}
I & 0 \\ 0 & C
\end{array} \right] \left[ \begin{array}{cc}
\dot{p} \\ \dot{L}
\end{array} \right] = \left[ \begin{array}{cc}
-R_p & 1 \\ -1 & -1/R_b
\end{array} \right] \left[ \begin{array}{cc}
\dot{p} \\ \dot{L}
\end{array} \right] + \left[ \begin{array}{cc}
0 \\ f_{in}
\end{array} \right] \ \ (\text{mode} \ 2)
\end{eqnarray*}
Supposing the system is initially in mode 1, the inflow causes the tank to start filling, which causes an outflow through resistance Rb. In this mode the outflow through the narrow pipe is zero. If L exceeds the level $L_{th}$, a switch from mode 1 to mode 2 occurs at the point in time when $L = L_{th}$.

\begin{figure}[h] \label{fig1} \vspace{-5mm}
\begin{center}
\includegraphics[scale=0.6]{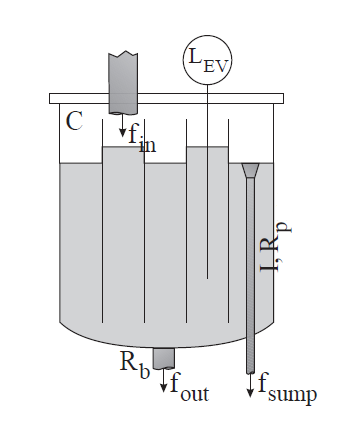}
\vspace{-3mm}
\caption{Schematic of the evaporator vessel}
\end{center} \vspace{-6mm}
\end{figure}

In the following, use the parameters $R_b = 1,R_p = 0.5, I = 0.5,C = 15, f_{in} = 0.25, L_{th} = 0.08$ and compute the following system matrices:
$$
\underline{\text{Mode 1}} : \bA_1 =  \left(\begin{array}{cc} -1 & 0\\ 0 & - \frac{1}{2} \end{array}\right), \ \ \bB_1 = \left(\begin{array}{c} 0\\ 1 \end{array}\right), \ \bC = \left(\begin{array}{cc} \frac{1}{2} & \frac{1}{2} \end{array}\right) 
$$
\vspace{-2mm}
$$
\underline{\text{Mode 2}}: \bA_2 =  \left(\begin{array}{cc} -1 & 2\\ - \frac{1}{2} & - \frac{1}{2} \end{array}\right)
, \ \ \bB_2 = \left(\begin{array}{c} 0\\ 1 \end{array}\right), \ \bC_2 = \left(\begin{array}{cc} \frac{1}{2} & \frac{1}{2} \end{array}\right) 
$$
First consider the following tuples of left and right interpolation nodes for each mode:
$$
\begin{cases} \blambda_1 = \{(-1.5),(-2,1)\} \\ \bmu_1 = \{(2),(0,0.5)\} \end{cases}, \ \ \begin{cases} \blambda_2 = \{(1),(1.5,-1.5)\} \\ \bmu_2 = \{(0),(2,-0.5)\} \end{cases}
$$
Hence, following the procedure described in Section 4, we recover the following system matrices:
$$
\underline{\text{Mode 1}}: \  \hat{\bE}_1 = \left(\begin{array}{cc} - \frac{1}{5} & - \frac{3}{20} \\[1mm] -1 & - \frac{1}{3} \end{array}\right), \ \hat{\bA}_1 =  \left(\begin{array}{cc} \frac{1}{10} & \frac{7}{60}\\[1mm] \frac{1}{2} & \frac{1}{6} \end{array}\right)	,  \ \hat{\bB}_1 = \left(\begin{array}{c} \frac{1}{5}\\ 1 \end{array}\right), \ \hat{\bC}_1 = \left(\begin{array}{cc} - \frac{1}{2} & - \frac{5}{12} \end{array}\right)
$$
$$
\underline{\text{Mode 2}}: \ \hat{\bE}_2 = \left(\begin{array}{cc} \frac{1}{2} & - \frac{5}{12}\\[1mm] \frac{1}{6} & - \frac{53}{360} \end{array}\right), \ \hat{\bA}_2 =  \left(\begin{array}{cc} - \frac{1}{2} & \frac{3}{8}\\[1mm] - \frac{4}{15} & \frac{17}{80} \end{array}\right),  \ \hat{\bB}_2 = \left(\begin{array}{c} 1\\[1mm] \frac{13}{30} \end{array}\right), \ \hat{\bC}_2 = \left(\begin{array}{cc} \frac{1}{2} & - \frac{3}{8} \end{array}\right)
$$
Note that the recovered realization is equivalent to the original one (no reduction has been enforced - the task was to recover the initial system). The coupling matrices are also computed:
\vspace{-2mm}
$$
\hat{\bK}_1 = \left(\begin{array}{cc} -1 & - \frac{1}{3} \\[1mm] - \frac{13}{30} & - \frac{17}{180} \end{array}\right), \ \ \hat{\bK}_2 = \left(\begin{array}{cc} \frac{11}{60} & - \frac{5}{36}\\[1mm] \frac{1}{2} & - \frac{5}{12} \end{array}\right)
$$

\vspace{-6mm}

\subsection{Second example}

For the next experiment, consider the CD player system from the SLICOT benchmark examples for MOR (see \cite{cvd02}). This linear system of order 120 has two inputs and two outputs. We consider that, at any given instance of time, only one input and one output are active (the others are not functional due to mechanical failure). For instance, consider mode $j$ to be activated whenever the $j^{th}$ input and the $j^{th}$ output are simultaneously failing (where $j \in \{1,2\})$.

In this way, we construct a LSS system with two operational modes. Both subsystems are stable SISO linear systems of order 120. This initial linear switched system $\Si$ will be 
reduced by means of the Loewner framework to obtain $\hat{\Si}_L$ and balanced truncation method proposed in \cite{mtc12} to obtain $\hat{\Si}_B$. 

The frequency response of each original subsystem is depicted in Fig.\;2.

\begin{figure}[h] \label{fig2} \vspace{-3mm}
\begin{center}
\includegraphics[scale=0.28]{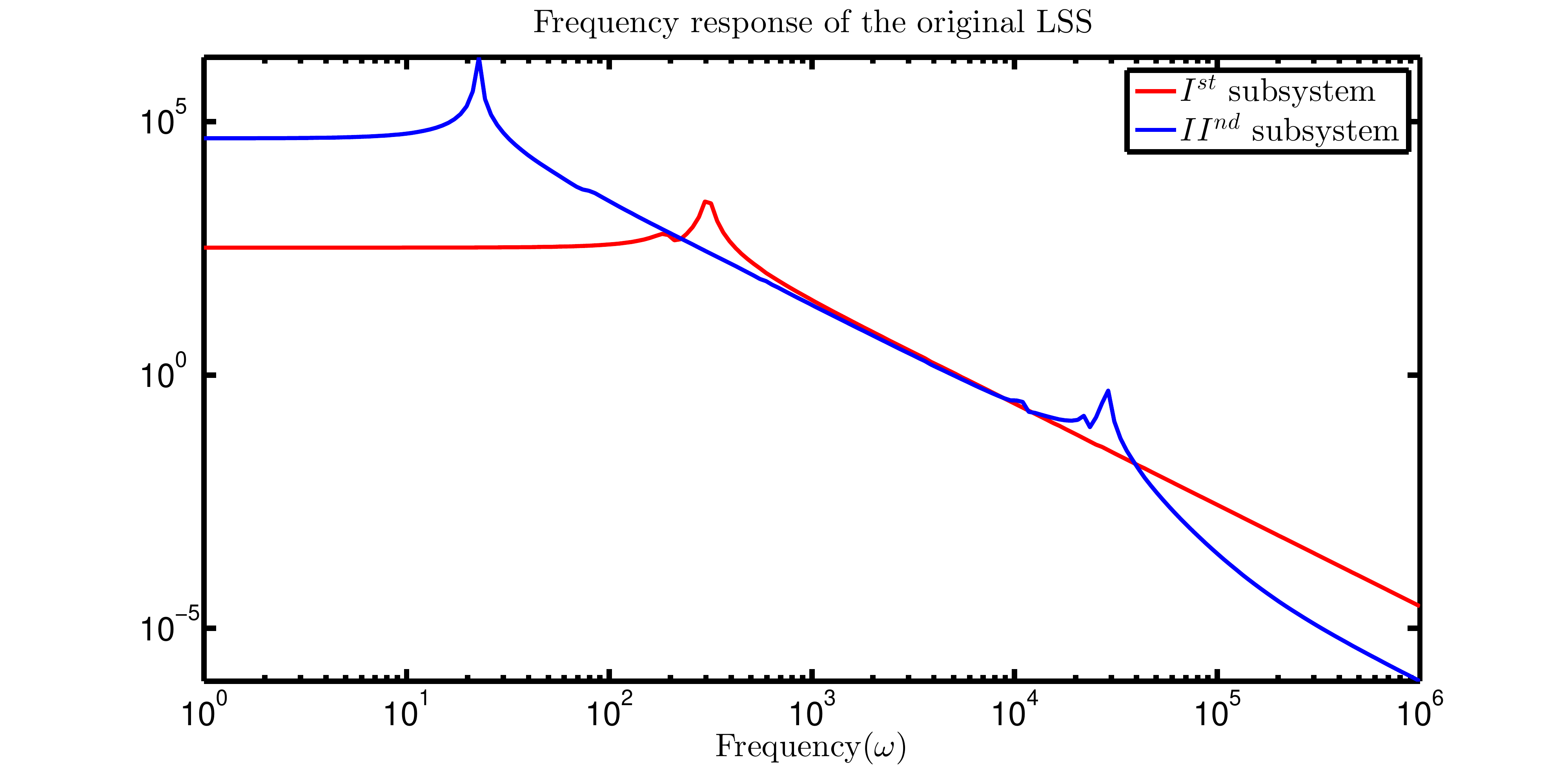}
\vspace{-3mm}
\caption{Frequency response of the original subsystems}
\end{center} \vspace{-6mm}
\end{figure}

For the Loewner method, we choose 160 logarithmically distributed  interpolation points inside $[10^1,10^5]j$. Fig.\;3 shows the decay of the singular values of the Loewner matrices corresponding to both subsystems. We notice that the $70^{th}$ singular values attain machine precision.

\begin{figure}[h] \label{fig3} \vspace{-3mm}
\begin{center}
\includegraphics[scale=0.28]{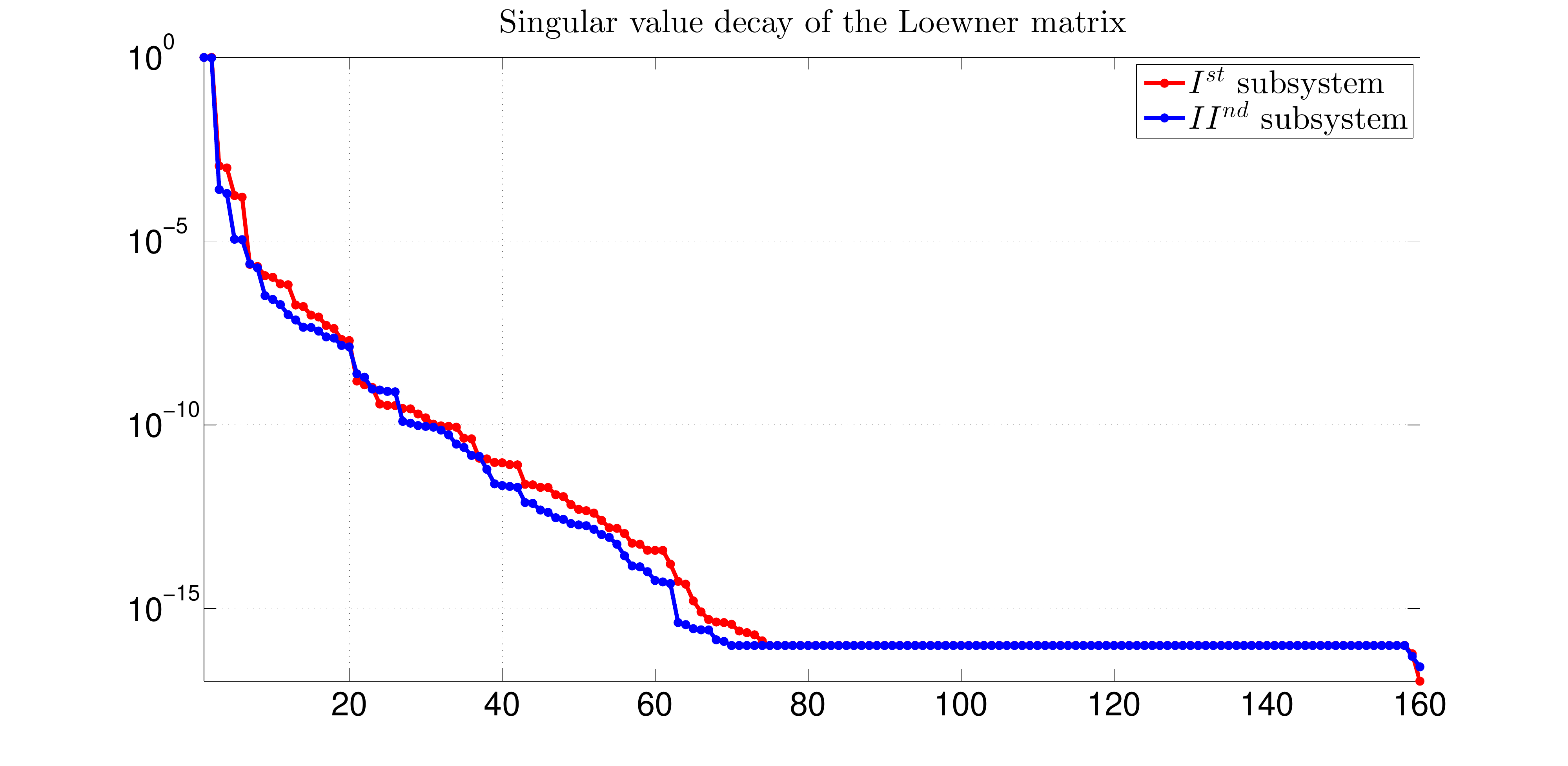}
\vspace{-6mm}
\caption{Decay of the singular values of the Loewner matrices}
\end{center} \vspace{-3mm}
\end{figure}
 
For $\hat{\Si}_L$ we decide to truncate at order $k = 28$ for both reduced systems. The same truncation order is chosen for $\hat{\Si}_B$. Next we compare the quality of approximation of the frequency response corresponding to the original system with the responses of the two reduced systems.

 In Fig.\;4 the relative error in frequency domain is depicted for both MOR methods (Loewner and BT). Notice that the Loewner method produces better results especially in the range of the selected sampling points. 

\begin{figure}[h] \label{fig4} \vspace{-1mm}
\begin{center}
\includegraphics[scale=0.3]{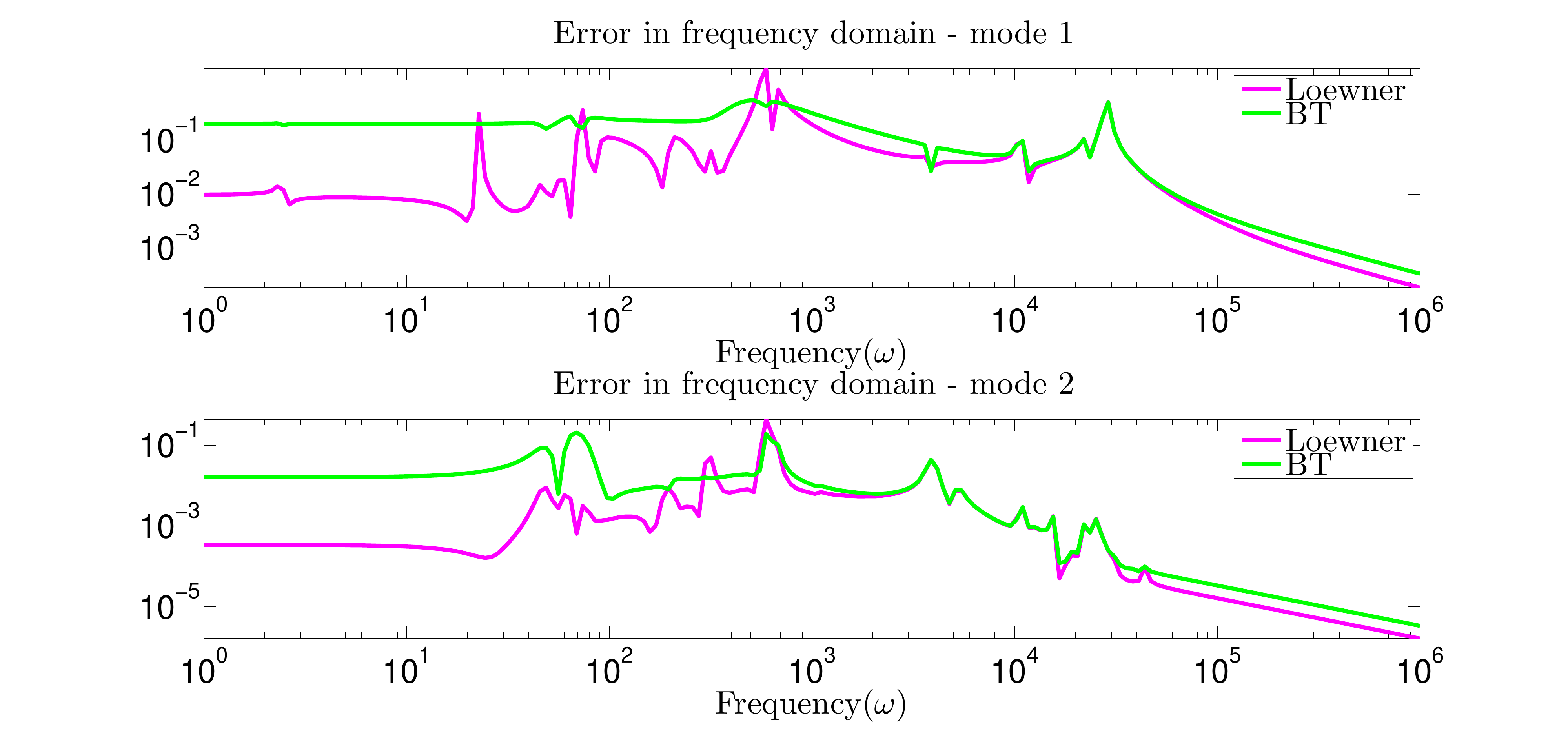}
\vspace{-3mm}
\caption{Frequency domain approximation error}
\end{center} \vspace{-3mm}
\end{figure}


Also, compare the time domain response of the original linear switched system against the ones corresponding to the two reduced models. We use a sinusoidal signal as the control input. The switching times are randomly chosen within [0,10]s. The blue rectangular signal in Fig.\;5 represents the switching signal. Notice that the output of the LSS is well approximated for both MOR methods, as it can be seen in the lower part of Fig.\;5.

\begin{figure}[h] \label{fig5} \vspace{-1mm}
\begin{center}
\includegraphics[scale=0.3]{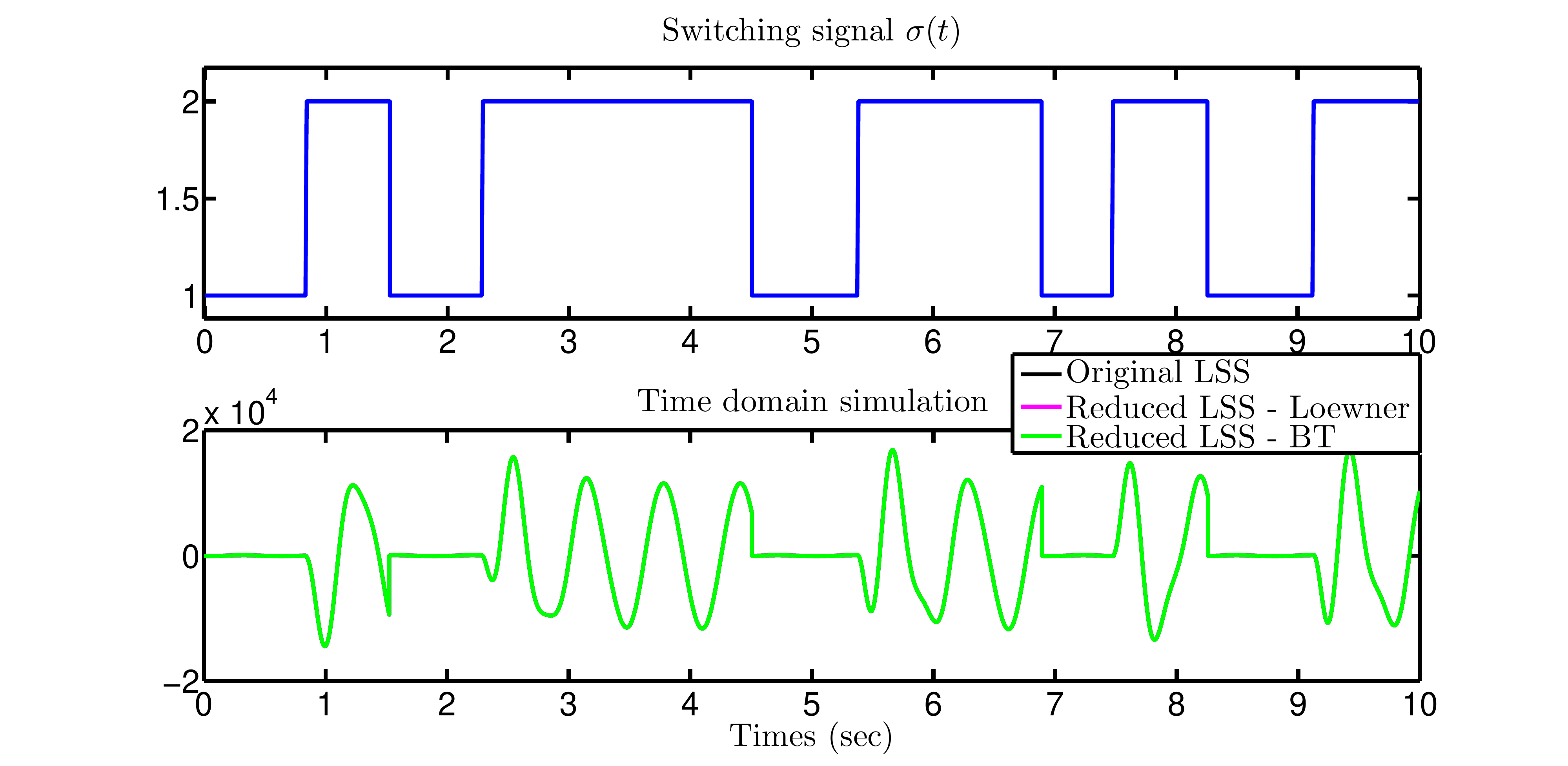}
\vspace{-3mm}
\caption{Time domain simulation}
\end{center} \vspace{-3mm}
\end{figure}

Finally, by inspecting the time domain error between the original response and the responses coming from the two reduced models (depicted in Fig.\;6), we notice that the Loewner method generally produces better results. The error curve corresponding to our proposed method is two orders of magnitude below the error curve corresponding to the BT method for most of the points on the time axis.

\begin{figure}[h] \label{fig6} \vspace{-1mm}
\begin{center}
\includegraphics[scale=0.26]{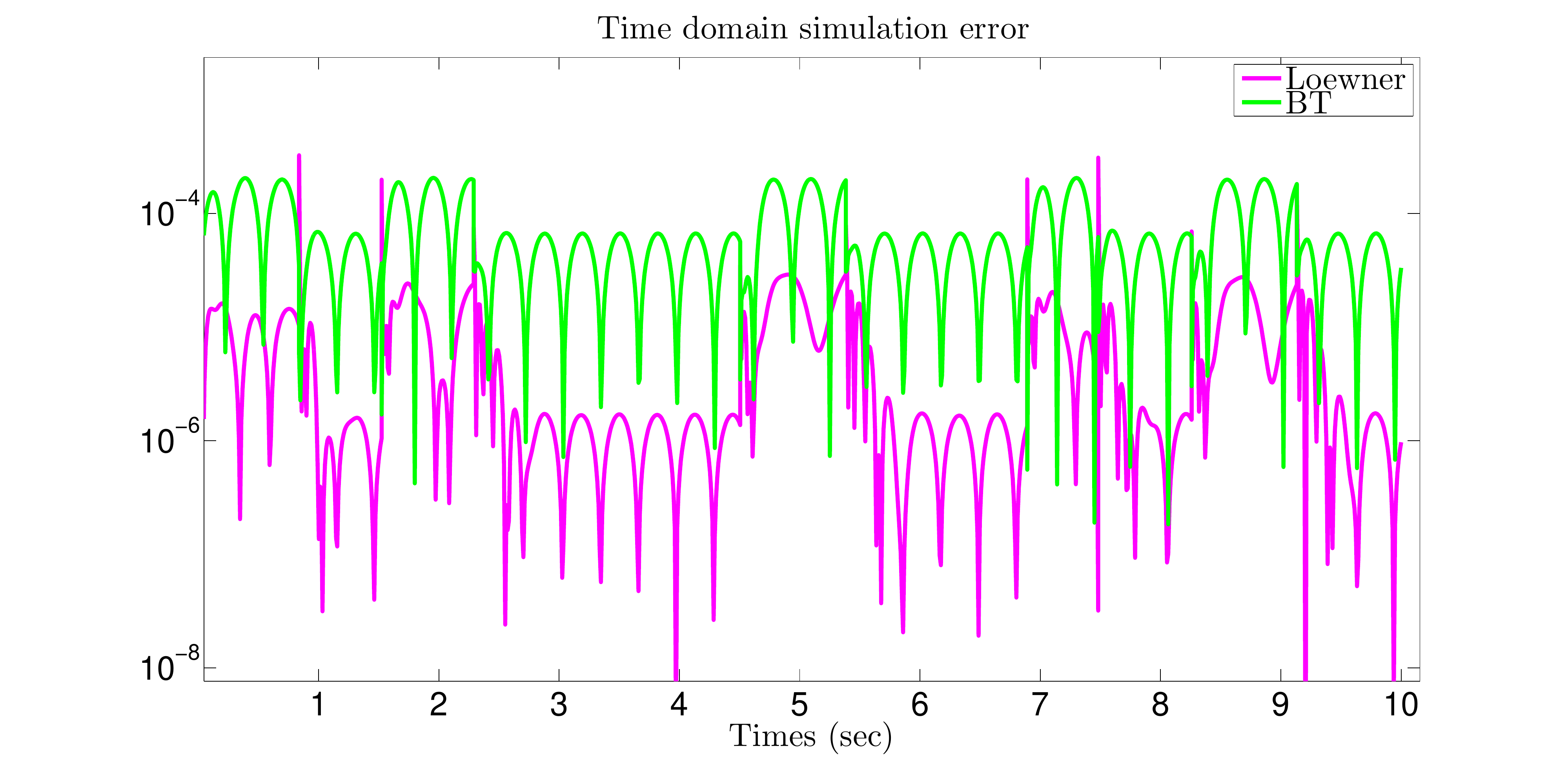}
\vspace{-3mm}
\caption{Time domain approximation error}
\end{center}
\end{figure} 

\subsection{Third example}

For the last experiment, consider a large scale LSS system constructed as in \cite{lsb14} from the original machine stand example given in \cite{ggm11}.

\vspace{-1mm}
\begin{figure}[h] \label{fig7} \vspace{-2mm}
\begin{center}
\includegraphics[scale=0.5]{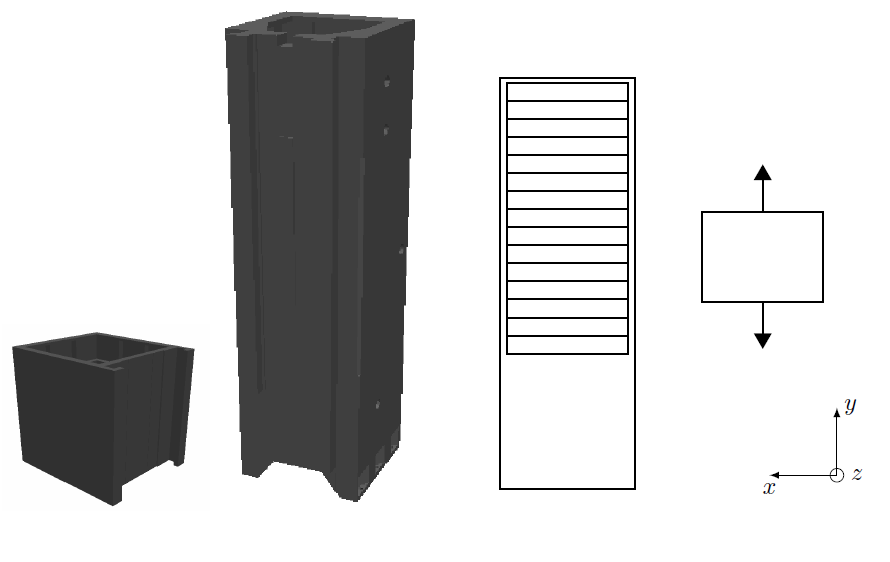}
\vspace{-3mm}
\caption{Schematic of the tool slide on the guide rails of the stand}
\end{center} \vspace{-3mm}
\end{figure}
\vspace{-1mm}

 In this example, the system variability is induced by a moving tool slide on the guide rails of the stand (see Fig.\;7). The aim is to determine the thermally driven displacement of the machine stand structure. Following the model setting in \cite{ggm11}, consider the heat equation with Robin boundary conditions. 
Using a finite element (FE) discretization and denoting
the external influences as the system input $z$, we obtain the dynamical heat model
\begin{equation} \label{eq:TV_model}
\bE_{th} \dot{\bx}(t) = \bA_{th}(t) \bx(t) + \bB_{th}(t) \bz(t)
\end{equation}
describing the deformation independent evolution of the
temperature field $\bx$ with the system matrices $\bE_{th},  \bA_{th}(t)$ and $ \bB_{th}(t)$. The variability of the model is described by time dependent matrices $\bA_{th}(t)$ and $\bB_{th}(t)$. This leads directly to the linear time varying system described by (\ref{eq:TV_model}). Since model reduction for LTV systems is a highly storage consuming procedure, the authors in \cite{lsb14} exploit properties of the spatially semi-discretized model to set up a LSS consisting of LTI subsystems only.
As described in \cite{ggm11}, the guide rails of
the machine stand are modeled as 15 equally distributed
horizontal segments (see Fig.\;7). Any of these segments is said to be completely covered by the tool slide if its midpoint lies within the height of the slide. On the other hand, each segment whose midpoint is not covered is treated as not in contact and therefore the slide always covers exactly 5 segments at each time. This in fact allows the stand to reach 11 distinct, discrete positions given by the model restrictions.
In this way, one can define the subsystems of the LSS as follows:
\begin{equation}
\Si_{\ell}: \begin{cases} \bE_{th} \dot{x} = \bA_{th}^{\ell} x+ \bB_{th}^{\ell} z^{\ell} \\ y = \bar{\bC} x, \end{cases}
\end{equation}
where $ \ell \in \{1, . . . , 11\}$. Note that the change
of the input operator Bth(t) is hidden in the input
zα itself, since it is sufficient to activate the correct
boundary parts by choosing the corresponding columns
in Bth via the input $z^\ell$. Therefore, the input operator
$B_{th}(t) := B_{th}$ becomes constant and the input variability is represented by the input $z^{\ell}$:
\begin{equation}
z_i^{\ell} : = \begin{cases} z_i, \ \text{segment i is in contact}, \\
0, \ \text{otherwise}, \end{cases}, \ \ i = 1, . . . , 15.
\vspace{-3mm}
\end{equation}
Here, $z_i \in \mathbb{R}$  is the thermal input as described in \cite{lsb14}. The only varying part influencing the model reduction process left in the dynamical system is the system ma
trix $\bA_{th}(t) := \bA^{\ell}_{th}$. 

After the finite element discretization was performed, we are given a SLS with 11 modes. Each subsystem has dimension $n=16626$. The $\bE$ and $\bC$ matrices are the same for all modes of the SLS. The $\bB$ matrices have 6 columns (corresponding to different inputs) and the C matrix has 9 rows (corresponding to different outputs).

For all the experiments performed, we take into consideration only two active modes (the first and the fifth). This corresponds to the particular case of $D=2$ covered in Section 4.
 
Also, although our new introduced method can be easily generalized to the MIMO case, we consider only (for simplicity reasons) the first input and the first output for each of the two modes  (the measurements used in the Loewner framework are hence scalar values).

Both subsystems are stable linear systems of order 16626 in sparse format. This initial large scale LSS will be again reduced (as in the second example) by means of the Loewner method and of the balanced truncation method proposed in \cite{mtc12}. The frequency response of each original subsystem is depicted in Fig.\;8.
 
\begin{figure}[h] \label{fig8} \vspace{-1mm}
\begin{center}
\includegraphics[scale=0.26]{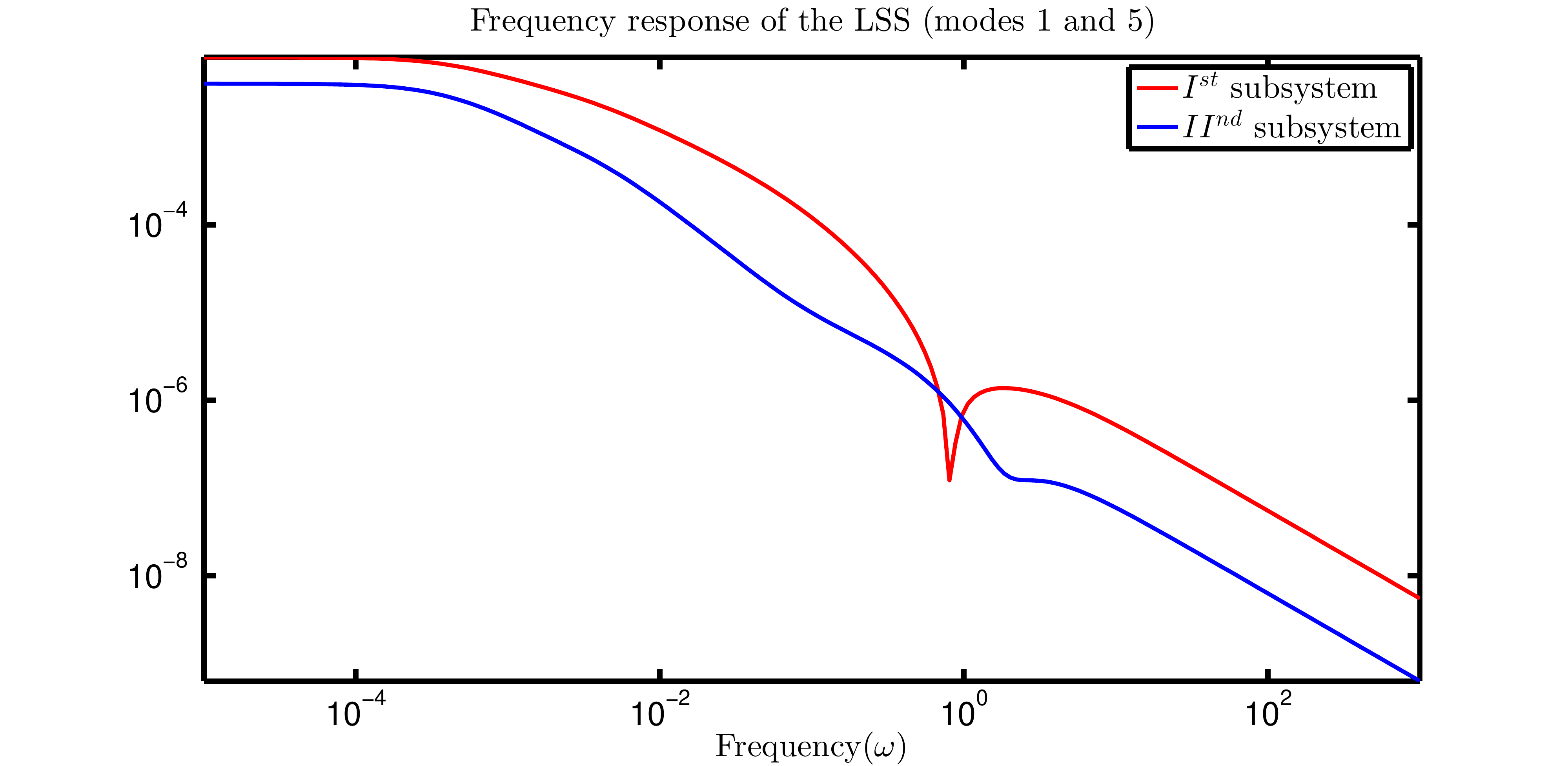}
\vspace{-3mm}
\caption{Frequency response of the original subsystems}
\end{center} \vspace{-3mm}
\end{figure}

For the Loewner method, we choose 200 logarithmically spaced interpolation points inside $[10^{-5},10^3]j$. The decay of the singular values of the Loewner matrices corresponding to both subsystems can be viewed in Fig.\;9. We notice that already the $70^{th}$ singular values attain machine precision.

\begin{figure}[h] \label{fig9} \vspace{-1mm}
\begin{center}
\includegraphics[scale=0.24]{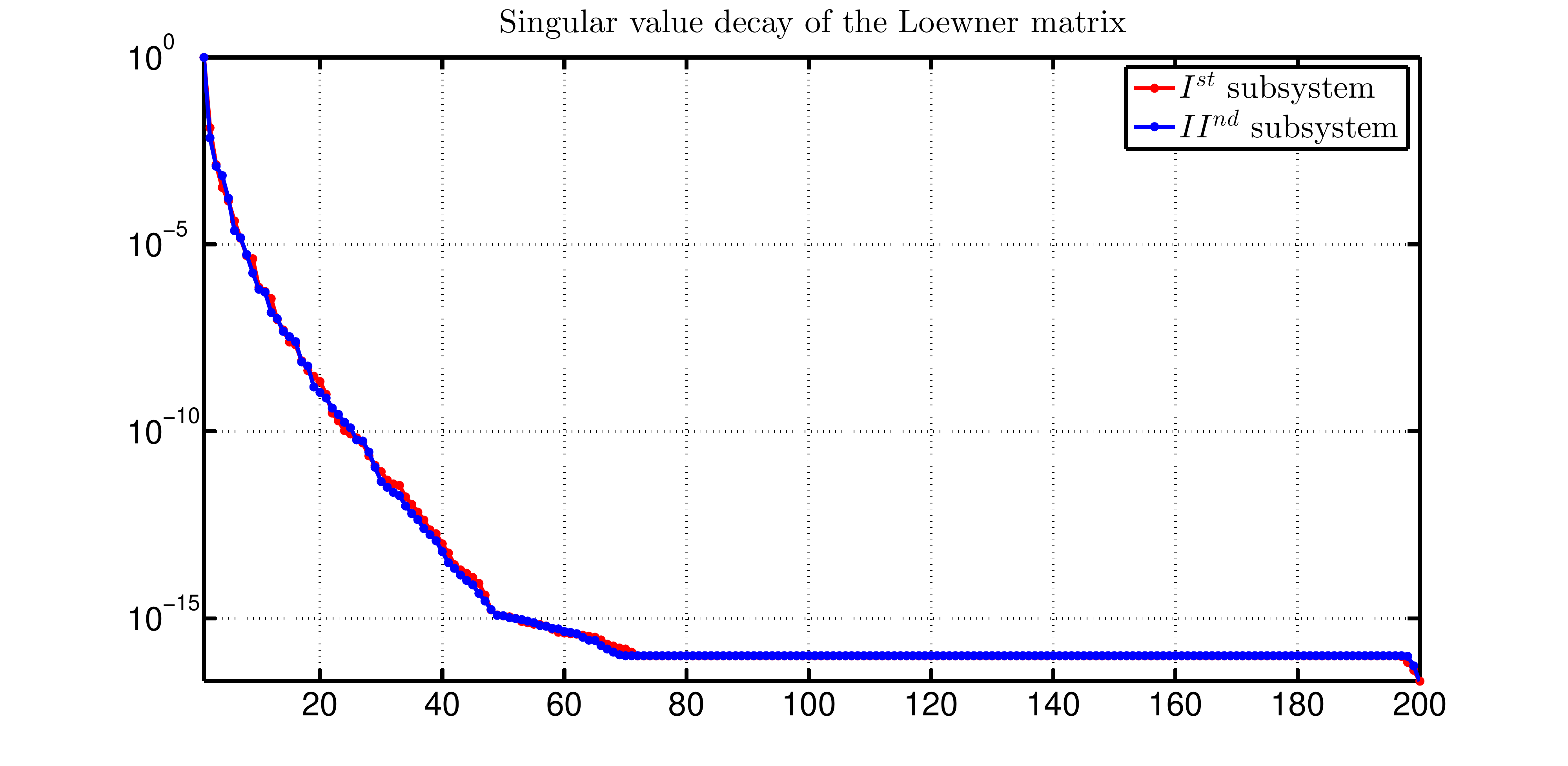}
\vspace{-5mm}
\caption{Decay of the singular values of the Loewner matrices}
\end{center} \vspace{-7mm}
\end{figure}

For the Loewner reduced order LSS (i.e $\Si_1$), we decide to truncate at order $k_1 = 26$ for both subsystems (which are stable LTI's). The same truncation order is chosen for the reduced order model computed via BT. Next we compare the quality of approximation of the frequency response. In Fig.\;10 the relative error in frequency domain is depicted for both MOR methods (Loewner and BT). Notice that the Loewner method generally produces better results.

\begin{figure}[h] \label{fig10} \vspace{-1mm}
\begin{center}
\includegraphics[scale=0.3]{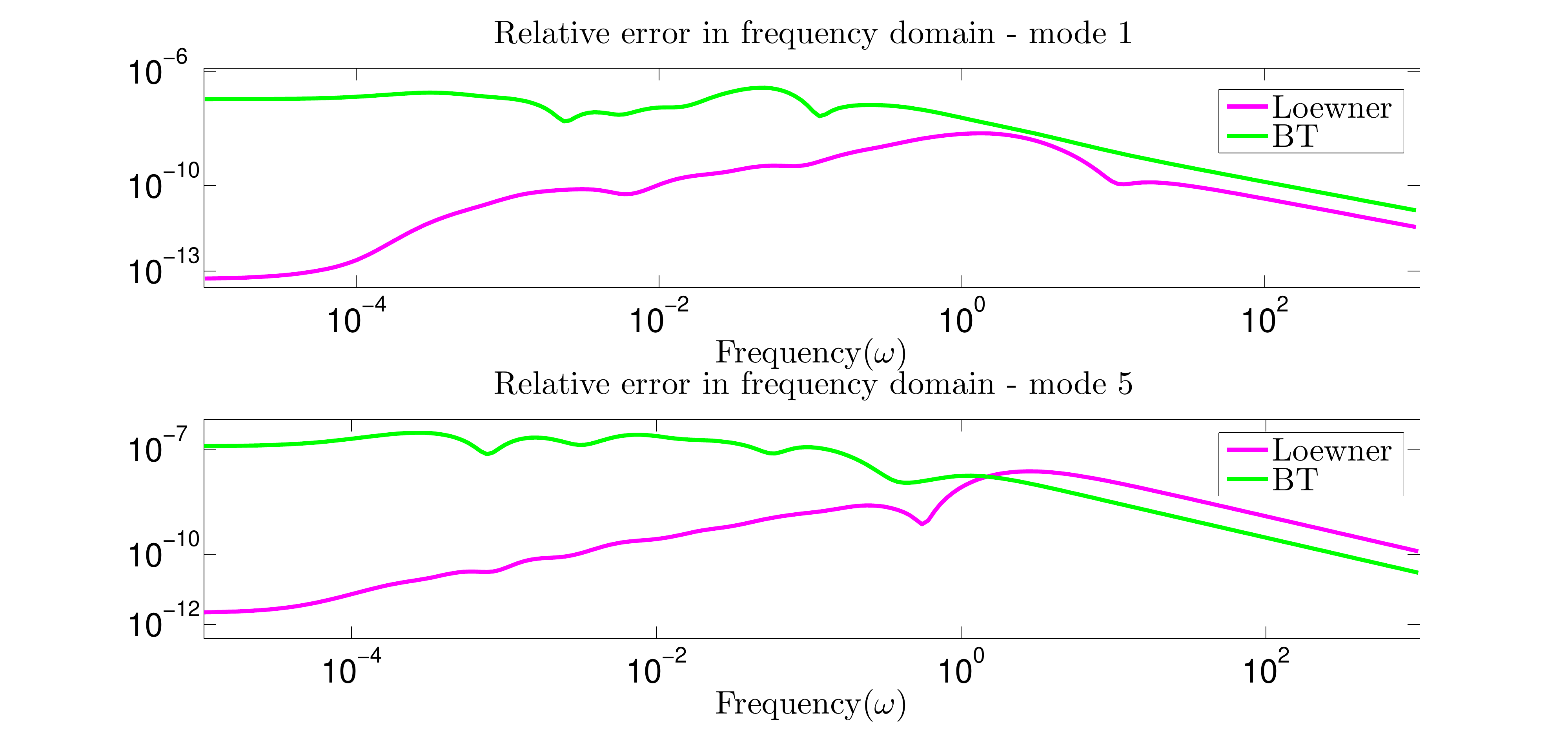}
\vspace{-3mm}
\caption{Frequency domain approximation error}
\end{center} \vspace{-7mm}
\end{figure}

Also, compare the time domain response of the original LSS against the ones corresponding to the two reduced models. The same configuration is used for the switching signal as in the second example. Notice that the output of the LSS is well approximated for both MOR methods, as it can be seen in the lower part of Fig.\;11.

\begin{figure}[h] \label{fig11} \vspace{-1mm}
\begin{center}
\includegraphics[scale=0.29]{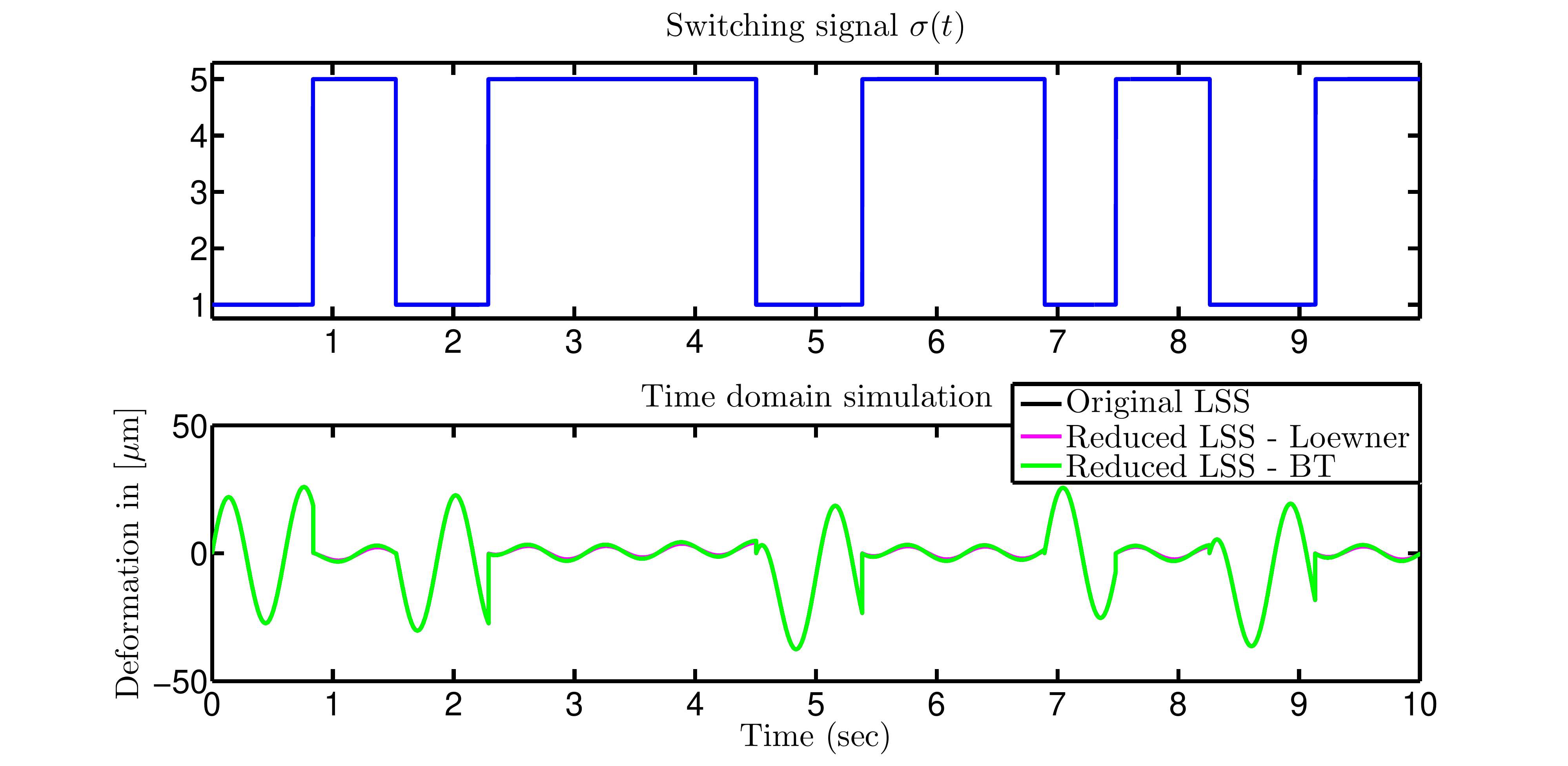}
\vspace{-3mm}
\caption{Time domain simulation} 
\end{center} \vspace{-7mm}
\end{figure}

Finally, by inspecting the time domain relative error between the original response and the responses coming from the two reduced models (depicted in Fig.\;12), we notice that the Loewner method generally produces better results. The error curve corresponding to our proposed method is below the error curve corresponding to the BT method for most of the points on the time axis.

\begin{figure}[h] \label{fig12} \vspace{-1mm}
\begin{center}
\includegraphics[scale=0.27]{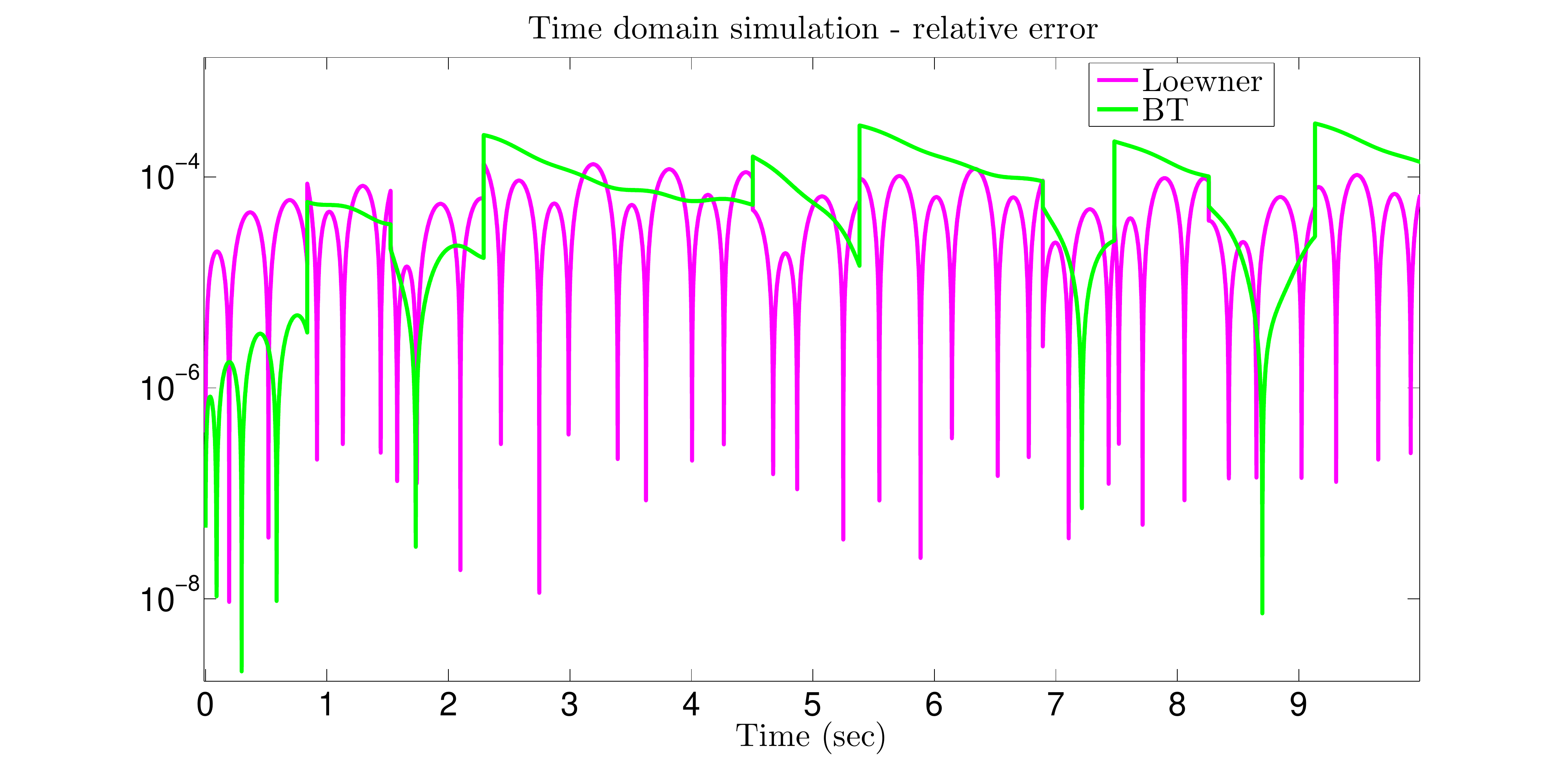}
\vspace{-3mm}
\caption{Time domain approximation error}
\end{center} \vspace{-8mm}
\end{figure}

\section{Summary and conclusions}

In this paper we address the problem of model reduction of linear switched systems from data consisting of values of high order transfer functions. The underlying philosophy of the Loewner framework is collect data and extract the desired information. Here he have extended this framework to the reduction of LSS. In general, for this type of systems, the data must be computed a priori, rather than measured (as for linear systems with no switching where one can use Vector Network Analyzers for instance). Having the required data, the next step would be to 
arrange it into matrix format. We have shown that the Loewner matrices (which basically represent the recovered E and A matrices of the underlying LSS) can be automatically calculated as solutions of Sylvester equations. In our framework, the transition/coupling matrices can be recovered from the given computed data as well. Since these matrices need not be square, they allow having different dimensions of the reduced state space in different modes.

 In a nutshell, given input/output data, we can construct with  no extensive computation, a singular high order model in generalized (descriptor) state space form. In applications the singular pencil $(\sIL,\IL)$ must be reduced at some stage. The {\it singular values} of the pencil $\,(\sIL,\,\IL)\,$ offer a {\it trade-off between accuracy of fit and complexity of the reduced system}. 
 
This approach to model reduction, first developed for linear time-invariant systems (see \cite{birkjour} for a survey), was later extended to {\it linear parametrized 
systems} \cite{ail12, IonAnt13, IonAnt14, iondiss} and to bilinear systems \cite{agi16}.

Three numerical examples demonstrate the effectiveness of the proposed approach. The quality of approximation for the reduced models was determined by performing both frequency and time domain tests. We have chosen a generalization of the classical balanced truncation method to LSS for comparison purposes. As opposed to most of the balancing methods we came across in the literature (\cite{ch09}, \cite{bgmb10}, \cite{sw11} and \cite{pwl13}), the method we choose (i.e \cite{mtc12}) does not require solving systems of LMI (linear matrix inequalities) which might be difficult for very large systems such as the one in Section 6.3. 
 The results of the new proposed method turned out to be better than the ones obtained when using the BT method.
\\

\noindent
\textbf{Acknowledgements}\\

The first and third authors would like to thank Mihaly Petreczky, Mert Bastug and Rafael Wisniewski for suggesting the problem of model reduction of LSS and for useful discussions. 

Also, the authors want to thank  Norman Lang for providing the LSS large scale system from Section 6.3 (used in \cite{lsb14}).

 Last but not least, we wish to thank Jens Saak for recommending the M-M.E.S.S. toolbox (see \cite{skb16}) for efficiently computing solutions of sparse large scale algebraic equations (which was very helpful for implementing the BT method for the $16626^{th}$ order system).

\bibliographystyle{siam}

\bibliography{GPA17_LSS}

\end{document}